\documentclass[12pt]{amsart}
\usepackage[utf8]{inputenc}
\usepackage[T1]{fontenc}
\usepackage{amsfonts}
\usepackage{amssymb}
\usepackage{amsmath}
\usepackage[english]{babel}
\usepackage{amsthm}
\theoremstyle{plain}
\newtheorem{theorem}{Theorem}[section]
\newtheorem{maintheorem}{Theorem}

\newtheorem{corollary}[theorem]{Corollary}

\newtheorem{lemma}[theorem]{Lemma}
\newtheorem{proposition}[theorem]{Proposition}
\theoremstyle{definition}

\newtheorem{remark}[theorem]{Remark}

\usepackage{geometry}
\numberwithin{equation}{section} 

\usepackage{lmodern,mathtools,microtype,needspace}

\usepackage[colorlinks, citecolor=blue,pagebackref,hypertexnames=false]{hyperref}

\newcommand{\R}{\mathbb R}
\newcommand{\C}{\mathbb C}
\newcommand{\one}{\mathbf1}

\newcommand{\Rea}{\mathfrak{Re}}

\newcommand{\cH}{\mathcal H}

\newcounter{comcount}

\DeclareMathOperator{\Ima}{\mathfrak{Im}}
\begin{document}
\title[Holomorphic calculus on $L^p(\ell^q)$]{Holomorphic calculus on $L^p(\ell^q)$ with optimal angle}
\author{Dangyang HE}
\date{}
\address{Dangyang He, Department of Mathematics, Sun Yat-sen (Zhongshan) University,\protect\newline Guangzhou, 510275, P.R. China}
\email{hedy28@mail.sysu.edu.cn}

\begin{abstract}
Let $\mathcal{L}$ generate a symmetric contraction semigroup. We prove a holomorphic H\"ormander calculus on $L^p(X;\ell^q)$ at the optimal angle $\max\{\arcsin|1-2/p|,\arcsin|1-2/q|\}$, with boundary Sobolev order $s>3$. Our approach combines a radial obstacle decomposition with contact inequalities and Mellin reconstruction, using one common contact set for all sequence coordinates and square-function scales. Sharpness is realized by a single generator built from finite Hamming cubes, with a translation embedding transferring the classical scalar obstruction to the fibre exponent.
\end{abstract}

\maketitle
\tableofcontents\enlargethispage{2pt}

\section{Introduction}\label{sec:introduction}
The spectral multiplier problem asks how regularity of a scalar
function on the spectrum controls boundedness of the corresponding
operator away from $L^2$. Let $(X,\mu)$ be a $\sigma$-finite measure
space, and let $\mathcal{L}$ be a nonnegative selfadjoint operator
on $L^2(X,\mu)$. Its spectral resolution $E_{\mathcal{L}}$ gives
\[
 \mathcal{L}=\int_0^\infty\lambda\,dE_{\mathcal{L}}(\lambda).
\]
For every bounded Borel function $m:[0,\infty)\to\mathbb C$,
the spectral theorem defines
\[
 m(\mathcal{L})=\int_0^\infty m(\lambda)\,dE_{\mathcal{L}}(\lambda),
 \qquad
 \|m(\mathcal{L})\|_{L^2\to L^2}
 \le \|m\|_\infty.
\]
For $p\ne2$, boundedness of $m(\mathcal{L})$ depends on both the
multiplier and the operator. Spectral multiplier theorems seek
conditions on $m$ which ensure its boundedness on $L^p(X)$.

We study this question for generators of symmetric contraction
semigroups. Throughout this paper, the consistent scalar actions
of $P_t=e^{-t\mathcal{L}}$ satisfy
\[
 \|P_t\|_{L^r(X)\to L^r(X)}\le1,
 \qquad t>0,\quad 1\le r\le\infty.
\]
The universal multiplier problem is to obtain estimates valid
for every such generator, with constants independent of the
generator and the underlying measure space. No metric or kernel
estimate enters into this formulation. Instead, holomorphy of
the multiplier supplies an analytic structure, and the opening
of its sector becomes a central parameter.

For $0<\theta<\pi$, put
\[
 \Sigma_\theta
 =\{z\in\mathbb C\setminus\{0\}:|\arg z|<\theta\},
 \qquad
 \mathbb C_+=\Sigma_{\pi/2}.
\]
Write $H^\infty(\Sigma_\theta)$ for the bounded holomorphic
functions on this sector, equipped with the supremum norm.
Its subspace $H^\infty_0(\Sigma_\theta)$ consists of functions
satisfying
\[
 |m(z)|\le C\min\{|z|^\varepsilon,|z|^{-\varepsilon}\}
\]
for some $C,\varepsilon>0$.
We prescribe the value $m(0)$ separately and write
$\Pi_0=\one_{\{0\}}(\mathcal{L})$, so that the contribution
from the kernel is $m(0)\Pi_0$. A bounded
$H^\infty(\Sigma_\theta)$ calculus on $L^p(X)$ means an estimate
\[
 \|m(\mathcal{L})\|_{L^p\to L^p}
 \le C_{p,\theta}
 \bigl(\|m\|_{H^\infty(\Sigma_\theta)}+|m(0)|\bigr)
\]
for every $m\in H^\infty(\Sigma_\theta)$ with its prescribed
value $m(0)$, where the constant is independent of $m$.

At a critical angle, bounded holomorphy is supplemented by
regularity on the boundary rays. Fix
$\eta\in C_c^\infty((0,\infty))$ equal to one near $[1/2,2]$.
Our Fourier convention and Sobolev norm are
\[
 \widehat h(\xi)=\int_{\mathbb R}e^{-ix\xi}h(x)\,dx,
 \qquad
 \|h\|_{H^s(\mathbb R)}^2
 =\frac1{2\pi}\int_{\mathbb R}
 (1+\xi^2)^s|\widehat h(\xi)|^2\,d\xi.
\]
For $s\ge0$, define
\begin{equation}
 \cH^s_\theta(m)
 =
 \|m\|_{H^\infty(\Sigma_\theta)}
 +
 \sum_{\epsilon=\pm1}\sup_{t>0}
 \|\eta(\cdot)m(t\,\cdot\,e^{i\epsilon\theta})\|_{H^s(\mathbb R)}.
\end{equation}
Boundary values are understood as nontangential traces, with
their continuous Sobolev representatives when $s>1/2$.
Different cutoffs of the stated kind give equivalent norms.
For Hilbert-valued functions and symbols, we use the same
definitions with the Hilbert norm.

Stein's semigroup theory \cite{Stein} established universal
multiplier estimates for symmetric diffusion semigroups,
including multipliers of Laplace transform type
\[
 m(\lambda)
 =\lambda\int_0^\infty a(t)e^{-t\lambda}\,dt,
 \qquad a\in L^\infty((0,\infty)).
\]
Cowling \cite{Cowling} proved that every symmetric contraction
generator has a bounded $H^\infty(\Sigma_\theta)$ calculus
on $L^p$ whenever
\[
 \theta>\pi\left|\frac12-\frac1p\right|,
 \qquad 1<p<\infty.
\]
This identifies an explicit sector using only contractivity.
The subsequent problem was to determine the smallest universal
angle and the regularity required on its boundary.

Imaginary powers provide a bridge between these two questions.
Here $\mathcal{L}^{i\tau}(I-\Pi_0)$ denotes the spectral
multiplier $\lambda^{i\tau}$ on $(0,\infty)$, with value zero
at the origin. Meda \cite{Meda1990} introduced a square-function
and Mellin-transform criterion which converts estimates for
these operators into multiplier theorems.
Its holomorphic formulations, developed by
Garc\'ia-Cuerva--Mauceri--Meda--Sj\"ogren--Torrea
\cite{GCMST} and Mauceri--Meda--Sj\"ogren \cite{MMS},
are recalled in \cite[Proposition~5]{CD}.
For example, if $0<\theta<\pi/2$ and $\beta\ge0$, an estimate
\[
 \|\mathcal{L}^{i\tau}(I-\Pi_0)\|_{L^p\to L^p}
 \le C(1+|\tau|)^\beta e^{\theta|\tau|}
\]
implies
\[
 \|m(\mathcal{L})\|_{L^p\to L^p}
 \le C\bigl(\cH^s_\theta(m)+|m(0)|\bigr),
 \qquad s>\beta+1.
\]
The exponential rate specifies the holomorphic sector, while
the polynomial growth controls the sufficient boundary
regularity. Cowling--Doust--McIntosh--Yagi \cite{CDMY}
developed the abstract relationship between square functions
and bounded holomorphic functional calculus.

For $1<r<\infty$, define
\[
 \sigma_r=\arcsin|1-2/r|,
 \qquad
 \phi_r=\frac{\pi}{2}-\sigma_r.
\]
Carbonaro--Dragi\v cevi\'c \cite{CD} determined the optimal
universal scalar calculus angle $\sigma_p$.
Their complex-time bilinear embedding theorem, proved through
a Bellman-function argument, yields
\[
 \|\mathcal{L}^{i\tau}(I-\Pi_0)\|_{L^p\to L^p}
 \le C_p(1+|\tau|)^{1/2}e^{\sigma_p|\tau|}.
\]
For $p\ne2$, their multiplier theorem takes the form
\[
 \|m(\mathcal{L})\|_{L^p\to L^p}
 \le C_{p,s}
 \bigl(\cH^s_{\sigma_p}(m)+|m(0)|\bigr),
 \qquad s>\frac32.
\]
The angle cannot be reduced uniformly over symmetric contraction
generators; the scalar obstruction was already known for the
Ornstein--Uhlenbeck operator \cite[Theorem~2]{GCMST}.
Kriegler--Weis \cite{KrieglerWeis} developed refined abstract
H\"ormander multiplier criteria through square-function
estimates. Haase--Pannasch \cite[Theorem~10.2]{HP}
extended this approach to general sectors and combined it with
the Carbonaro--Dragi\v cevi\'c estimate to obtain the scalar
boundary order $s>1$ at the same angle $\sigma_p$.

A complementary real-variable theory allows multipliers which
need not be holomorphic. Its regularity conditions involve
localized Sobolev quantities such as
\[
 \sup_{t>0}
 \|\eta(\cdot)m(t\,\cdot)\|_{H^s(\mathbb R)}.
\]
H\"ormander's Fourier multiplier theorem
\cite{Hormander1960} is a fundamental result in this direction.
Duong--Ouhabaz--Sikora \cite{DOS2002} connected sharp spectral
multiplier estimates with Gaussian heat-kernel bounds and
Plancherel-type estimates.
Sikora--Yan--Yao \cite{SYY2014} developed this framework under
generalized Gaussian estimates and Plancherel or restriction-type
conditions. These results relate the required regularity to
geometric and spectral information about the operator.
Our universal holomorphic estimates use only symmetric
contractivity.

For vector-valued functions, the geometry of the fibre enters
the multiplier problem. We consider the concrete mixed spaces
$L^p(X;\ell^q)$, $1<p,q<\infty$, with norm
\[
 \|f\|_{L^p(X;\ell^q)}
 =
 \left(
 \int_X
 \left(\sum_j|f_j(x)|^q\right)^{p/q}
 \,d\mu(x)
 \right)^{1/p}.
\]
For a Banach space $E$ and $1\le p<\infty$, we also use
\[
 \|f\|_{L^{p,\infty}(X;E)}
 =
 \sup_{\lambda>0}
 \lambda\,\mu\{x:\|f(x)\|_E>\lambda\}^{1/p}.
\]
Thus weak norms refer to the distribution of the fibre norm
in the $X$ variable. The real semigroup acts coordinatewise
and contractively on these mixed spaces. This follows from
the positive Dunford--Schwartz modulus of each scalar operator,
as recalled in Subsection~\ref{subsec:radial}.

The spaces $\ell^q$, $1<q<\infty$, are UMD Banach lattices.
Hyt\"onen \cite{Hytonen2007} developed
Littlewood--Paley--Stein theory for symmetric diffusion
semigroups on UMD-valued $L^p$ spaces.
Xu \cite{Xu} proved maximal inequalities for UMD lattices
and obtained a bounded holomorphic calculus at some angle
below $\pi/2$ for suitable complex interpolation spaces
between Hilbert spaces and UMD spaces; these include all
UMD Banach lattices. These classical results provide the
initial nonsharp calculus and vector-valued estimates used
below. Our common-contact obstacle argument yields the
optimal upper angle, and the cube construction gives the
matching obstruction.
In his later work, Xu \cite[Theorems~1.2 and~1.4]{Xu2025}
developed a holomorphic-calculus approach to Banach-valued
Littlewood--Paley--Stein inequalities for Poisson and analytic
semigroups, with quantitative estimates governed by
martingale type and cotype.
Under generalized Gaussian estimates,
Deleaval--Kemppainen--Kriegler \cite{DKK2021}
established real H\"ormander multiplier theorems on UMD
lattice-valued spaces.

For the Ornstein--Uhlenbeck family, van Neerven
\cite{vanNeerven2026} obtained a dimension-free
$H^\infty$ calculus at some angle below $\pi/2$ for
arbitrary UMD fibres. Independently, van Neerven obtained the optimal angle of
dimension-free analyticity for $L^q(0,1)$-valued
Ornstein--Uhlenbeck semigroups for the outer exponent $p=2$
in an unpublished manuscript dated 11 August 2026
\cite{vanNeervenOUAugust}.
In the same setting, his subsequent manuscript dated
2 September 2026 \cite{vanNeervenOUSeptember}
determines the optimal angle for every $1<p<\infty$.
The present paper establishes imaginary-power and holomorphic
multiplier estimates for general symmetric contraction
generators through a radial obstacle decomposition,
common-contact inequalities and Mellin reconstruction.

Our main result shows that the analyticity angle determined by van Neerven for the Ornstein--Uhlenbeck family \cite{vanNeervenOUSeptember} is also optimal uniformly over all symmetric contraction generators acting on $L^p(X;\ell^q)$.
We establish imaginary-power and holomorphic multiplier estimates at the corresponding calculus angle. We write
\begin{equation*}
\sigma_{p,q}=\max\{\sigma_p,\sigma_q\},
 \qquad
 \phi_{p,q}
 =\frac{\pi}{2}-\sigma_{p,q}
 =\min\{\phi_p,\phi_q\},
\end{equation*}
where $\sigma_{p,q}$ is the optimal universal calculus angle and $\phi_{p,q}$ is the corresponding optimal universal angle of semigroup analyticity. Thus both the outer exponent and the fibre exponent can determine the sector.

The proof uses an obstacle decomposition with one common contact set. Ouyang--Spector--Stockdale \cite{OSS2026} introduced an obstacle decomposition for the fractional Laplacian in their study of the vector Riesz transform. Chen--Jiang--Li--Li \cite{CJLL} developed an abstract semigroup version for positive selfadjoint operators with sub-Markovian semigroups. In these decompositions, the exceptional set is the support of a potential, and the bounded good part saturates the obstacle there. We construct a radial version for $L^p$ sources with values in uniformly smooth reflexive Banach
spaces. Applied in $\ell^q(L^2((0,\infty),ds/s))$, it produces one contact set shared by all sequence coordinates and square-function scales. The construction and its limiting argument are proved in detail in Section~\ref{sec:contact}.

We record the remaining conventions before stating the theorem.
For $1<r<\infty$, write $r'=r/(r-1)$.
Complex pairings are linear in the first argument.
For Bochner functions, an unsubscripted pairing denotes
the integrated duality
\[
 \langle F,G\rangle
 =
 \int_X\langle F(x),G(x)\rangle_E\,d\mu(x),
 \qquad
 F\in L^p(X;E),\quad
 G\in L^{p'}(X;E^*).
\]
We abbreviate $L^r(X;E)$ to $L^r(E)$ when $X$ is fixed.
Functional calculi initially agree with the scalar spectral
calculus on finite sequences with entries in $L^2$;
their extensions are understood by consistency and density.

\begin{maintheorem}\label{thm:A}
Let $\mathcal{L}$ be as above and let $1<q<\infty$.
For $1<p<\infty$ and $\tau\in\R$,
\begin{equation}\label{eq:power-main}
 \|\mathcal{L}^{i\tau}(I-\Pi_0)\|_{L^p(X;\ell^q)\to L^p(X;\ell^q)}
 \le C_{p,q}(1+|\tau|)^2e^{\sigma_{p,q}|\tau|}.
\end{equation}
If $\sigma_{p,q}>0$, $s>3$, and
$m\in H^\infty(\Sigma_{\sigma_{p,q}})$ has
$\cH^s_{\sigma_{p,q}}(m)<\infty$, then
\begin{equation}\label{eq:horm-main}
 \|m(\mathcal{L})\|_{L^p(X;\ell^q)\to L^p(X;\ell^q)}
 \le C_{p,q,s}
 \bigl(\cH^s_{\sigma_{p,q}}(m)+|m(0)|\bigr).
\end{equation}
Consequently the generator has a bounded
$H^\infty(\Sigma_\theta)$ calculus for every
$\theta>\sigma_{p,q}$.
At $p=q=2$, every bounded Borel multiplier is allowed.
Moreover, the angle $\sigma_{p,q}$ is optimal uniformly
over symmetric contraction generators.
\end{maintheorem}

The influence of the fibre is already visible when $p=2$:
the optimal universal angle on $L^2(X;\ell^4)$ is $\pi/6$,
although the scalar $L^2$ calculus is supplied directly by
the spectral theorem. Thus the fibre can impose a positive
angular threshold even when the outer exponent is Hilbertian.
On the diagonal $p=q$, Fubini reduces the result to the scalar
theory; for $q=2$, Gaussian amplification does the same.
All sequence estimates have constants independent of the
number of coordinates.

Section~\ref{sec:applications} records two consequences of
this angular information. Corollary~\ref{cor:rotated-square}
gives two-sided lattice square-function estimates on
complex-time rays, uniformly on every closed subsector
of $\Sigma_{\phi_{p,q}}$.
Corollary~\ref{cor:sectorial-maximal} gives the corresponding
sectorial maximal estimate for every
$0<\gamma<\phi_{p,q}$, with the optimal universal angular
threshold. These applications combine the classical
square-function and maximal arguments with the full strict
mixed-norm angular range supplied by Theorem~\ref{thm:A}.

The result applies, in particular, to Laplace--Beltrami
operators and symmetric Markov generators. It also gives
explicit mixed-norm estimates for fractional and higher-order
diffusion. For $\beta>0$, let $\mathcal{L}^\beta$ denote
the spectral power. If $\beta\sigma_{p,q}<\pi/2$, then,
for every $0<\gamma<\pi/2-\beta\sigma_{p,q}$,
\[
 \sup_{z\in\Sigma_\gamma}
 \|e^{-z\mathcal{L}^\beta}\|_
 {L^p(X;\ell^q)\to L^p(X;\ell^q)}
 \le C_{p,q,\beta,\gamma}.
\]
To see this, choose $\sigma_{p,q}<\theta<\pi$ such that
$\gamma+\beta\theta<\pi/2$ and apply the calculus to
$m_z(\lambda)=e^{-z\lambda^\beta}$, using the principal branch
of $\lambda^\beta$ and prescribing $m_z(0)=1$.
Its $H^\infty(\Sigma_\theta)$ norm is at most one, uniformly
for $z\in\Sigma_\gamma$.
The permitted powers include $\beta>1$ when $p$ and $q$
are sufficiently close to $2$.

The boundary order $s>3$ and the polynomial factor in
\eqref{eq:power-main} are not asserted to be optimal.
In particular, our result does not improve the scalar
regularity thresholds described above.

\medskip
\noindent\emph{Strategy of the proof.}
The main step is to construct a radial obstacle in
$\ell^q(L^2(ds/s))$, with one common contact set for all
sequence coordinates and square-function scales.
This common set allows the contact inequalities to be
integrated directly in the mixed norm, with constants
independent of the number of coordinates and of the
scale interval.

Mellin reconstruction converts these contact estimates into
weak bounds for imaginary powers, while square-function
analysis and synthesis remove the auxiliary scale variable.
The contact estimate is independent of the specific multiplier,
which enters through the bounds for its Mellin transform.
Interpolation between nearby outer exponents then preserves
the optimal exponential growth rate.
The boundary multiplier theorem follows from these
imaginary-power estimates and the abstract holomorphic
H\"ormander calculus.

For sharpness, we construct one symmetric contraction generator,
fixed independently of $p$ and $q$, as a direct sum of finite
Hamming-cube generators. The scalar obstruction is classical
\cite{IvanisviliNazarov}. A translation embedding transfers the
obstruction at exponent $q$ to the fibre of every
$L^p(X;\ell^q)$, while a single-coordinate test detects the
obstruction at exponent $p$. Thus the same generator accounts
for both contributions to the optimal angle.
Moreover, no finite increase in boundary Sobolev regularity
permits a universally smaller angle.

\medskip
\noindent\emph{Organization of the paper.}
Section~\ref{sec:contact} contains the auxiliary arguments:
Subsection~\ref{subsec:radial} constructs the radial obstacle,
Subsection~\ref{subsec:complex} proves complex contractivity,
Subsection~\ref{subsec:contact} establishes the common energy
and contact estimates, and Subsection~\ref{subsec:mellin}
contains the Mellin reconstruction and mixed-norm integration
lemmas.
In Section~\ref{sec:calculus},
Subsection~\ref{subsec:upper} proves the imaginary-power and
multiplier estimates, and Subsection~\ref{sec:sharpness}
establishes sharpness.
Section~\ref{sec:applications} derives the square-function
and sectorial maximal consequences.

\section{Preparation}\label{sec:contact}

\subsection{A radial obstacle decomposition}\label{subsec:radial}

We retain the notation of Section~\ref{sec:introduction}. Recall that a Dunford--Schwartz operator is a complex-linear operator
whose consistent actions on $L^1(X)$ and $L^\infty(X)$ are contractions.
Its intermediate $L^r$ actions are the consistent contractions supplied
by interpolation.  We shall use the classical linear-modulus fact:
every such operator $T$ is dominated by a positive Dunford--Schwartz
operator $S$, so that $|Tv|\le S|v|$ for $v\in L^1(X)$.
This is the single-operator assertion in
\cite[Lemma~3.4]{Taggart}. For a separable Banach space $E$ and a finite-rank simple function $v$, domination implies
\[
 \|Tv(x)\|_E\le S(\|v\|_E)(x).
\]
Indeed, for every $e^*\in E^*$ of norm at most one,
\[
 |\langle Tv,e^*\rangle|
 =|T\langle v,e^*\rangle|
 \le S(|\langle v,e^*\rangle|)
 \le S(\|v\|_E).
\]
Take the supremum over a countable norming set for the separable range of $v$. Approximation by finite-rank simple functions proves the assertion for $v\in L^r(X;E)$, $1\le r<\infty$, and defines the consistent Banach-valued extensions. In particular, for a positive $S$ and $v=(v_j)$,
\begin{equation}
 \|(Sv_j(x))_j\|_{\ell^q}\le S(\|v\|_{\ell^q})(x),
\end{equation}
and
\[
 \|Tv\|_{L^r(X;E)}\le\|v\|_{L^r(X;E)}.
\]
For $v\in L^1(X;E)\cap L^\infty(X;E)$, the same pointwise domination shows
\[
 \|Tv\|_{L^\infty(E)}\le\|v\|_{L^\infty(E)}.
\]
This intersection action is sufficient for every $L^\infty$ bound below. 

\begin{lemma}\label{lem:erg-threshold}
For a normed space $E$ and $c>0$, put
\begin{equation*}
S_cz=\left(1-\frac c{\|z\|_E}\right)_+z,
 \qquad S_c0=0.
\end{equation*}
Then $\|z\|_E=\|S_cz\|_E+\|z-S_cz\|_E$ and $\|z-S_cz\|_E=\min\{\|z\|_E,c\}$. Moreover, $S_c$ is $1$-Lipschitz:
\begin{equation*}
\|S_cz-S_cw\|_E\le\|z-w\|_E.
\end{equation*}

\end{lemma}

\begin{proof}
The two norm identities follow from the definition. Suppose without loss of generality that $\|z\|_E \ge \|w\|_E$. If $\|z\|_E\le c$, then by definition $S_cz=S_cw=0$. If $\|w\|_E \le c<\|z\|_E$, then it is clear that
\begin{equation*}
    \| S_cz-S_cw\|_E =\|z\|_E-c\le \|z\|_E-\|w\|_E \le \|z-w\|_E.
\end{equation*}
If instead $c<\|w\|_E\le \|z\|_E$, a direct computation yields $S_cz-S_cw = (1-c/\|z\|_E)(z-w)+c(1/\|w\|_E-1/\|z\|_E)w$. The triangle inequality then gives
\begin{align*}
\| S_cz-S_cw\|_E &\le \left( 1-\frac{c}{\|z\|_E} \right) \|z-w\|_E + \frac{c}{\|z\|_E} \left(\|z\|_E-\|w\|_E\right)\\
&\le \left( 1-\frac{c}{\|z\|_E} \right) \|z-w\|_E + \frac{c}{\|z\|_E}  \|z-w\|_E\\
&\le \|z-w\|_E.\qedhere
\end{align*}

\end{proof}
Recall that $E$ is uniformly smooth if
\[
 \lim_{t\downarrow0}\frac1t
 \sup_{\|v\|_E=\|z\|_E=1}
 \left(\frac{\|v+tz\|_E+\|v-tz\|_E}{2}-1\right)=0.
\]
We use its classical consequences for norming functionals and Bochner-space duality maps; see \cite[Chapter~I, Sections~3--4, and Chapter~II, Section~4]{Cioranescu}.

The obstacle method of Ouyang--Spector--Stockdale
\cite{OSS2026} arises from a fractional-Laplacian obstacle problem in their study of the vector Riesz transform. Chen--Jiang--Li--Li \cite{CJLL} developed an abstract decomposition for positive selfadjoint operators with sub-Markovian semigroups.
The next theorem adapts this mechanism to $L^p$ sources with values in a uniformly smooth reflexive Banach space.

\begin{theorem}\label{lem:radial}
Let $E$ be uniformly smooth and reflexive, let
$1<p<r<\infty$, and consider the spaces $L^s(X;E)$.
Suppose $-\mathcal{A}$ generates consistent strongly continuous
semigroups on $L^p(X;E)$ and $L^r(X;E)$, with
\[
 \|e^{-t\mathcal{A}}\|_{L^s(X;E)\to L^s(X;E)}
 \le e^{-\kappa t},
 \qquad s=p,r,
\]
for some $\kappa>0$. Assume also that for every convex
increasing function
$\Psi:[0,\infty)\to[0,\infty)$ with $\Psi(0)=0$,
\[
 \int_X\Psi(\|e^{-t\mathcal{A}}u(x)\|_E)\,d\mu
 \le
 \int_X\Psi(\|u(x)\|_E)\,d\mu,
\]
whenever $u\in L^p(X;E)\cap L^r(X;E)$ and the right-hand
side is finite.

For $f\in L^p(X;E)\cap L^r(X;E)$ and $h>0$, there exist $u\in\operatorname{Dom}_{L^p}(\mathcal{A})
\cap\operatorname{Dom}_{L^r}(\mathcal{A})$ and $g\in L^p(X;E)\cap L^r(X;E)$ such that
\[
 \mathcal{A}u+g=f,\qquad
 \|g(x)\|_E\le h,\qquad
 g(x)=h\frac{u(x)}{\|u(x)\|_E}
 \quad\hbox{where }u(x)\ne0,
\]
and
\[
\begin{aligned}
 \|g\|_{L^p(X;E)}&\le\|f\|_{L^p(X;E)},\\
 \mu\{u\ne0\}&\le h^{-p}\|f\|_{L^p(X;E)}^p,\\
 \|g\|_{L^r(X;E)}^r
 &\le h^{r-p}\|f\|_{L^p(X;E)}^p.
\end{aligned}
\]
\end{theorem}

\begin{proof}
\noindent\emph{Step 1. The decomposition of $f$.}
Write $\sigma_E(v)=v/\|v\|_E$ for $v\ne0$, with
$\sigma_E(0)=0$. Define
\[
 S_c(v)=\left(1-\frac{c}{\|v\|_E}\right)_+v,
 \qquad S_c(0)=0.
\]
By Lemma~\ref{lem:erg-threshold}, this map is
$1$-Lipschitz. The Laplace resolvent formula \cite[Theorem~II.1.10(ii)--(iii)]{EngelNagel} and the semigroup bound give
\[
 \|(\mathcal{A}+d)^{-1}\|_{L^s(X;E)\to L^s(X;E)}
 \le(\kappa+d)^{-1},
 \qquad s=p,r,\quad d>0.
\]
For $\varepsilon>0$, put $d=\varepsilon^{-1}$ and
consider the mapping
\[
 \Phi_\varepsilon(v)
 :=(\mathcal{A}+\varepsilon^{-1})^{-1}
       \bigl(f+\varepsilon^{-1}S_{\varepsilon h}(v)\bigr).
\]
For $s=p,r$, Lemma~\ref{lem:erg-threshold} guarantees $\|S_{\varepsilon h}(v)\|_{L^s(X;E)}
 \le\|v\|_{L^s(X;E)}$. Thus $\Phi_\varepsilon$ maps $L^s(X;E)$ into itself and
\[
 \|\Phi_\varepsilon(v)-\Phi_\varepsilon(w)\|_{L^s(X;E)}
 \le
 \frac{\varepsilon^{-1}}{\varepsilon^{-1}+\kappa}
 \|v-w\|_{L^s(X;E)}.
\]
Apply the Banach contraction mapping theorem on
$L^p(X;E)\cap L^r(X;E)$, equipped with the sum norm.
Consistency of the resolvents makes the mapping
well-defined on this intersection. We obtain a unique
fixed point $u_\varepsilon$ satisfying
\[
 u_\varepsilon
 =(\mathcal{A}+\varepsilon^{-1})^{-1}
       \bigl(f+\varepsilon^{-1}
                    S_{\varepsilon h}(u_\varepsilon)\bigr),
\]
and $\|u_\varepsilon\|_{L^s(X;E)}
 \le\kappa^{-1}\|f\|_{L^s(X;E)}$ for $s=r,p$. Since the resolvent takes values in the operator domain, $u_\varepsilon\in\operatorname{Dom}_{L^p}(\mathcal{A})
\cap\operatorname{Dom}_{L^r}(\mathcal{A})$. Setting
\[
 g_\varepsilon
 =\varepsilon^{-1}
       \bigl(u_\varepsilon-S_{\varepsilon h}u_\varepsilon\bigr),
\]
we obtain $\mathcal{A}u_\varepsilon+g_\varepsilon=f$. Moreover, $\|g_\varepsilon(x)\|_E\le h$ by Lemma~\ref{lem:erg-threshold}.

\medskip
\noindent\emph{Step 2. The entropy test.}
For $v\ne0$, denote by $j_E(v)$ its norming functional,
normalized by
\[
 \|j_E(v)\|_{E^*}=1,\qquad
 \langle v,j_E(v)\rangle=\|v\|_E,
\]
and set $j_E(0)=0$. Uniform smoothness makes this
functional unique. Put
\[
 J_{p,E}(v)=\|v\|_E^{p-1}j_E(v),\qquad J_{p,E}(0)=0.
\]
Thus $\|J_{p,E}(v)\|_{E^*}=\|v\|_E^{p-1}$, and $\langle v,J_{p,E}(v)\rangle=\|v\|_E^p$.

Next, consider the convex function $\Psi_\varepsilon(t)
 =\int_0^t\min(s/\varepsilon,h)^{p-1}\,ds$. Equivalently,
\[
 \Psi_\varepsilon(t)=
 \begin{cases}
 \displaystyle\frac{t^p}{p\varepsilon^{p-1}},
       &0\le t\le\varepsilon h,\\[6pt]
 \displaystyle h^{p-1}t-\frac{p-1}{p}\varepsilon h^p,
       &t\ge\varepsilon h.
 \end{cases}
\]
One checks straightforwardly
\[
 \Psi_\varepsilon(t)
 \le\frac{t^p}{p\varepsilon^{p-1}},
 \qquad
 \Psi_\varepsilon'(t)
 =\min(t/\varepsilon,h)^{p-1}.
\]
Observe that
\[
 g_\varepsilon
 =\min\bigl(\|u_\varepsilon\|_E/\varepsilon,h\bigr)
       \sigma_E(u_\varepsilon).
\]
Thus $g_\varepsilon$ is a positive scalar multiple of
$u_\varepsilon$ wherever $u_\varepsilon\ne0$, and hence
$j_E(g_\varepsilon)=j_E(u_\varepsilon)$ there. It follows
that
\begin{equation}\label{eq:lem:erg:1}
 \Psi_\varepsilon'(\|u_\varepsilon\|_E)
       j_E(u_\varepsilon)
 =J_{p,E}(g_\varepsilon).
\end{equation}
Both sides vanish where $u_\varepsilon=0$.

Convexity of $\Psi_\varepsilon$ and the norming property
give
\[
\begin{split}
 \Psi_\varepsilon(\|z\|_E)
 &\ge
 \Psi_\varepsilon(\|u\|_E)
 +\Psi_\varepsilon'(\|u\|_E)
       \bigl(\|z\|_E-\|u\|_E\bigr)\\
 &\ge
 \Psi_\varepsilon(\|u\|_E)
 +\Psi_\varepsilon'(\|u\|_E)
       \Rea\langle z-u,j_E(u)\rangle.
\end{split}
\]
Replacing $u$ by $u_\varepsilon$, integrating and using
\eqref{eq:lem:erg:1}, we obtain
\[
\begin{split}
 \Rea\int_X
 \langle z-u_\varepsilon,J_{p,E}(g_\varepsilon)\rangle
 \,d\mu \le
 \int_X\Psi_\varepsilon(\|z\|_E)\,d\mu
 -\int_X\Psi_\varepsilon(\|u_\varepsilon\|_E)\,d\mu.
\end{split}
\]
Take $z=e^{-t\mathcal{A}}u_\varepsilon$. The assumption makes the
right-hand side non-positive, so
\[
 \Rea\int_X
 \left\langle
 \frac{u_\varepsilon-e^{-t\mathcal{A}}u_\varepsilon}{t},
 J_{p,E}(g_\varepsilon)
 \right\rangle\,d\mu
 \ge0.
\]
Notice that $\|g_\varepsilon(x)\|_E
 \le\varepsilon^{-1}\|u_\varepsilon(x)\|_E$.
Thus $g_\varepsilon\in L^p(X;E)$ and
$J_{p,E}(g_\varepsilon)\in L^{p'}(X;E^*)$.
Since $u_\varepsilon\in\operatorname{Dom}_{L^p}(\mathcal{A})$, one deduces $\frac{u_\varepsilon-e^{-t\mathcal{A}}u_\varepsilon}{t}
 \longrightarrow \mathcal{A}u_\varepsilon$ in $L^p(X;E)$ as $t\downarrow0$. H\"older's inequality therefore permits passage to the limit $t\downarrow0$ and gives
\[
 \Rea\int_X
 \langle \mathcal{A}u_\varepsilon,J_{p,E}(g_\varepsilon)\rangle
 \,d\mu
 \ge0.
\]

\medskip
\noindent\emph{Step 3. Uniform estimates for $g_\varepsilon$.}
Since $f=\mathcal{A}u_\varepsilon+g_\varepsilon$, the preceding
inequality yields
\[
\begin{split}
 \|g_\varepsilon\|_{L^p(X;E)}^p
 &=\Rea\int_X
       \langle g_\varepsilon,J_{p,E}(g_\varepsilon)\rangle
       \,d\mu\\
 &=\Rea\int_X
       \langle f,J_{p,E}(g_\varepsilon)\rangle\,d\mu
   -\Rea\int_X
       \langle \mathcal{A}u_\varepsilon,J_{p,E}(g_\varepsilon)\rangle
       \,d\mu\\
 &\le
 \Rea\int_X
       \langle f,J_{p,E}(g_\varepsilon)\rangle\,d\mu\\
 &\le
 \|f\|_{L^p(X;E)}
 \|g_\varepsilon\|_{L^p(X;E)}^{p-1}.
\end{split}
\]
Consequently, $\|g_\varepsilon\|_{L^p(X;E)}\le\|f\|_{L^p(X;E)}$.
Moreover, since $r>p$ and
\[
 \|g_\varepsilon(x)\|_E\le h,
\]
one has
\[
 \|g_\varepsilon\|_{L^r(X;E)}^r
 \le h^{r-p}\|g_\varepsilon\|_{L^p(X;E)}^p
 \le h^{r-p}\|f\|_{L^p(X;E)}^p.
\]

\medskip
\noindent\emph{Step 4. Strong convergence of $u_\varepsilon$.}
We now let $\varepsilon\downarrow0$. In this step, write $\|\cdot\|_s$ for the norm of $L^s(X;E)$ and let $J_{p,E}$ act pointwise on functions.

We use the classical smoothness and duality-map properties of Bochner spaces; see \cite[Chapter~I, Sections~3--4, and Chapter~II, Section~4]{Cioranescu}. The pointwise map $J_{p,E}:L^p(X;E)\to L^{p'}(X;E^*)$ is uniformly continuous on bounded sets. In particular, the functional $v\mapsto\|v\|_p^p/p$ is continuously Fr\'echet differentiable over the underlying real space, with
\begin{equation}\label{eq:Bochner}
\lim_{\substack{t\to0\\t\in\mathbb R}}
 \frac{\|v+tz\|_p^p-\|v\|_p^p}{pt}
 =
 \Rea\int_X\langle z,J_{p,E}(v)\rangle\,d\mu.
\end{equation}
Set $w_\varepsilon
 =u_\varepsilon-\varepsilon g_\varepsilon
 =S_{\varepsilon h}u_\varepsilon$.
The radial formulas show that $g_\varepsilon
 =h w_\varepsilon/\|w_\varepsilon\|_E$ if $w_\varepsilon \ne0$. On $\{w_\varepsilon=0\}$, one has $\|g_\varepsilon\|_E\le h$. It follows that, for every $\lambda>0$, $w_\varepsilon
 =S_{\lambda h}(w_\varepsilon+\lambda g_\varepsilon)$. Therefore, Lemma~\ref{lem:erg-threshold} ensures
\[
 \|w_\varepsilon-w_\delta\|_p
 \le
 \|w_\varepsilon-w_\delta
       +\lambda(g_\varepsilon-g_\delta)\|_p.
\]
Raising both sides to the power $p$, subtracting $\|w_\varepsilon-w_\delta\|_p^p$ from both sides, dividing by $p\lambda$ and letting $\lambda \downarrow 0$, one obtains by \eqref{eq:Bochner},
\begin{equation*}
\Rea\int_X\langle g_\varepsilon-g_\delta, J_{p,E}(w_\varepsilon-w_\delta)\rangle\,d\mu \ge0.
\end{equation*}
Next, for $v\in\operatorname{Dom}_{L^p}(\mathcal{A})$, contractivity means that
$\|e^{-t\mathcal{A}}v\|_p^p\le e^{-p\kappa t}\|v\|_p^p$. Hence
\[
 \frac{\|e^{-t\mathcal{A}}v\|_p^p-\|v\|_p^p}{pt}
 \le \frac{e^{-p\kappa t}-1}{pt}\|v\|_p^p.
\]
Since $\frac{e^{-t\mathcal{A}}v-v}{t}\longrightarrow -\mathcal{A}v$ in $L^p(X;E)$, \eqref{eq:Bochner} yields
\[
 \Rea\int_X\langle \mathcal{A}v,J_{p,E}(v)\rangle\,d\mu
 \ge\kappa\|v\|_p^p.
\]
Apply this inequality to $v=u_\varepsilon-u_\delta$.
Since $\mathcal{A}(u_\varepsilon-u_\delta)=-(g_\varepsilon-g_\delta)$, combining the preceding inequalities shows
\begin{align*}
\kappa\|u_\varepsilon-u_\delta\|_p^p
 &\le
 -\Rea\int_X
 \langle g_\varepsilon-g_\delta,
         J_{p,E}(u_\varepsilon-u_\delta)\rangle\,d\mu\\
 &\le
 \Rea\int_X
 \left\langle g_\varepsilon-g_\delta,\,
 J_{p,E}(w_\varepsilon-w_\delta)
 -J_{p,E}(u_\varepsilon-u_\delta)\right\rangle\,d\mu\\
 &\le
 2\|f\|_p
 \left\|
 J_{p,E}(w_\varepsilon-w_\delta)
 -J_{p,E}(u_\varepsilon-u_\delta)
 \right\|_{L^{p'}(X;E^*)}.
\end{align*}
For $0<\varepsilon,\delta\le1$, it is plain that $\|(w_\varepsilon-w_\delta)
       -(u_\varepsilon-u_\delta)\|_p\le
 \varepsilon\|g_\varepsilon\|_p
 +\delta\|g_\delta\|_p\le(\varepsilon+\delta)\|f\|_p
 \to 0$. Both $w_\varepsilon-w_\delta$ and $u_\varepsilon-u_\delta$ lie in the ball of radius $2\kappa^{-1}\|f\|_p$: indeed, $\|w_\varepsilon\|_p\le\|u_\varepsilon\|_p$ by radial thresholding. Uniform continuity of the duality map on this ball therefore shows $\|J_{p,E}(w_\varepsilon-w_\delta)-J_{p,E}(u_\varepsilon-u_\delta)\|_{L^{p'}(X;E^*)} \to0$. Thus $u_\varepsilon$ is Cauchy and converges to a function $u$ in $L^p(X;E)$.

\medskip
\noindent\emph{Step 5. Regularity.}
Steps~1 and~3 bound both $u_\varepsilon$ and
$g_\varepsilon$ uniformly in $L^p(X;E)\cap L^r(X;E)$.
This intersection, equipped with the sum norm, is
reflexive. The weak subsequence theorem for reflexive Banach spaces \cite[Theorem~3.18]{Brezis} therefore applies to the pairs $(u_\varepsilon,g_\varepsilon)$. Hence, along a sequence $\varepsilon_n\downarrow0$, $u_{\varepsilon_n}\to u$, $g_{\varepsilon_n}\to g$ weakly in both $L^p(X;E)$ and $L^r(X;E)$. For $s=p,r$, the generator $\mathcal{A}$ on $L^s(X;E)$ is closed \cite[Theorem~II.1.4]{EngelNagel}. Its graph is a norm-closed linear subspace, hence is weakly closed by \cite[Theorem~3.7]{Brezis}. Thus passing weakly to the limit in $\mathcal{A}u_{\varepsilon_n}=f-g_{\varepsilon_n}$ shows that $(u,\,f-g)$ belongs to the $L^s(X;E)$ graph of $\mathcal{A}$ and therefore
\[
 u\in\operatorname{Dom}_{L^p}(\mathcal{A})
       \cap\operatorname{Dom}_{L^r}(\mathcal{A}),
 \qquad \mathcal{A}u+g=f.
\]
The subset $\{v\in L^p(X;E):\|v(x)\|_E\le h\ \hbox{a.e.}\}$ is norm closed and convex, hence weakly closed by the same theorem. It follows that $\|g(x)\|_E\le h$ almost everywhere. Weak lower semicontinuity also guarantees $\|g\|_p \le \|f\|_p$.

It remains to identify $g$ on the support of $u$.
Passing to a subsequence if necessary, let $u_{\varepsilon_n}(x)\to u(x)$ in $E$ almost everywhere. Let $\Omega=\{u\ne0\}$. For a.e. $x\in\Omega$, one has, for $n$ large enough,
\[
 \|u_{\varepsilon_n}(x)\|_E
 \ge\frac12\|u(x)\|_E
 >\varepsilon_n h.
\]
For such $n$,
\[
 g_{\varepsilon_n}(x) = \min\left(\frac{\|u_{\varepsilon_n}(x)\|_E}{\varepsilon_n},\, h\right) \sigma_E(u_{\varepsilon_n})
 =h\frac{u_{\varepsilon_n}(x)}
          {\|u_{\varepsilon_n}(x)\|_E}
 \to h\frac{u(x)}{\|u(x)\|_E}.
\]
For every measurable $F\subset\Omega$ with
$\mu(F)<\infty$, the bound $\|g_{\varepsilon_n}-h u/\|u\|_E\|_E\le2h$ on $F$ and dominated convergence imply $\mathbf1_F g_{\varepsilon_n} \to \mathbf1_F h\frac{u}{\|u\|_E}$ strongly in $L^p(X;E)$. Meanwhile, the same sequence converges weakly to $\mathbf1_F g$. Multiplication by $\mathbf1_F$ is bounded and linear, so it preserves weak convergence. Uniqueness of weak limits identifies these limits. Since $X$ is sigma-finite, exhaust $\Omega$ by such sets $F$ to obtain
\[
 g=h\frac{u}{\|u\|_E}
\]
a.e. on $\Omega$. Finally, since $\|g(x)\|_E=h$ a.e. on $\Omega$, 
\begin{align*}
    h^p \mu(\Omega) = \int_\Omega \|g(x)\|_E^p \,d\mu \le \|g\|_p^p\le \|f\|_p^p.
\end{align*}
Also, for $r>p$, the pointwise upper bound implies
\[
 \|g\|_r^r=\int_X\|g(x)\|_E^{r-p}\|g(x)\|_E^p\,d\mu
 \le h^{r-p}\|g\|_p^p\le h^{r-p}\|f\|_p^p.
\]
This completes the proof of Theorem~\ref{lem:radial}.
\end{proof}

\begin{remark}
The positive constant $\kappa$ is used both in the
fixed-point construction and in the strong convergence of the
approximate potentials. The resulting potential satisfies
\[
 \|u\|_{L^s(X;E)}
 \le \kappa^{-1}\|f\|_{L^s(X;E)},
 \qquad s=p,r,
\]
whereas the estimates for the good part and the exceptional
set are independent of $\kappa$. This distinction is useful in
Proposition~\ref{prop:common-contact} below. For a compact interval $I\subset(0,\infty)$, the constant there is
\[
 \kappa_I=\inf_{s\in I}s^{-1/2}>0,
\]
which tends to zero along an increasing exhaustion of
$(0,\infty)$. Nevertheless, the constants in the good-part,
exceptional-set, energy and contact estimates are independent
of $I$. These uniform estimates permit the later passage to
all scales.
\end{remark}

\subsection{The complex contractivity}\label{subsec:complex}
Put
\[
 c_r=\frac{|r-2|}{2\sqrt{r-1}}=\tan\sigma_r.
\]
Throughout this subsection, we write $c_{p,q}=\max\{c_p,c_q\}$.

The sharp mixed-duality-map inequality in
Lemma~\ref{lem:umd-mixed-sector}, including its optimal constant and the differential mean and covariance calculation, also appears in
\cite[Lemma~4.1]{vanNeervenOptimalOU}. For Ornstein--Uhlenbeck semigroups, van Neerven combines this inequality with the symmetric Mehler coupling to
obtain complex contractivity. We retain the calculation in our notation for completeness.
Haase's general two-point reduction \cite[Theorem~2.2]{HaaseForms} then allows us to treat arbitrary symmetric contraction generators.
The resulting Hilbert-valued extension will be used in the common-contact argument below.

\begin{lemma}\label{lem:umd-mixed-sector}
Let $1<p,q<\infty$ and $N\ge1$. Define for $v\in \mathbb C^N$
\begin{equation*}
J_{p,q}(v)_j=\|v\|_q^{p-q}|v_j|^{q-2}v_j,
\end{equation*}
with $|v_j|^{q-2}v_j=0$ when $v_j=0$ and $J_{p,q}(0)=0$. Then, for $u,v\in\mathbb C^N$,
\begin{equation}\label{eq:umd-two-point-sector}
\begin{aligned}
 &\left|\Ima\sum_j(u_j-v_j)
       \left(\overline{J_{p,q}(u)_j-J_{p,q}(v)_j}\right) \right|\\
 &\qquad\le c_{p,q}\Rea\sum_j(u_j-v_j)
       \left(\overline{J_{p,q}(u)_j-J_{p,q}(v)_j}\right).
\end{aligned}
\end{equation}
The real part on the right is nonnegative. The constant $c_{p,q}$ is
optimal when $N\ge2$; for $N=1$ the optimal constant is $c_p$.
\end{lemma}

\medskip
\noindent\emph{Remark.}
For $E=\ell_N^q$, the map $J_{p,q}$ is precisely the duality
map $J_{p,E}$ defined in the proof of Theorem~\ref{lem:radial}.
The two subscripts merely display the outer and fibre exponents.

\begin{proof}
First suppose that $z_j\ne0$ for every $j$, put $r=\|z\|_q$, and write $h_j/z_j=a_j+ib_j$. Introduce the probability weights $\mu_j=\frac{|z_j|^q}{r^q}$. These weights satisfy $\sum_j \mu_j=1$. The expectation and variance of a real vector $\omega\in \mathbb R^N$ are
\begin{equation*}
    \mathbb E_\mu (\omega) = \sum_j \mu_j \omega_j,\qquad \mathrm{Var}_\mu (\omega) = \sum_j \mu_j \left( \omega_j-\mathbb E_\mu(\omega)\right)^2=\mathbb E_\mu(\omega^2) - \mathbb E_\mu(\omega)^2.
\end{equation*}
Also, the covariance of $\omega, \nu\in \mathbb R^N$ is
\begin{equation*}
    \mathrm{Cov}_\mu(\omega,\nu) = \mathbb E_\mu \left( (\omega-\mathbb E_\mu(\omega)) (\nu-\mathbb E_\mu(\nu)) \right) = \mathbb E_\mu (\omega \nu) - \mathbb E_\mu(\omega) \mathbb E_\mu(\nu).
\end{equation*}
It follows by Cauchy-Schwarz that $|\mathrm{Cov}_\mu(\omega,\nu)|\le \mathrm{Var}_\mu (\omega)^{1/2} \mathrm{Var}_\mu (\nu)^{1/2}$.

Differentiating with respect to the real parameter $t$ gives
\begin{equation*}
\begin{aligned}
 (DJ_{p,q}(z)h)_j
 &:= \frac{d}{dt} J_{p,q}(z+th)_j \bigg|_{t=0}\\
 &= r^{p-q}|z_j|^{q-2}z_j
       \left[ (q-1)a_j+ib_j+(p-q)\mathbb E_\mu a \right].
\end{aligned}
\end{equation*}
With $\mathcal C(z,h)=\sum_jh_j\overline{(DJ_{p,q}(z)h)_j}$, it follows that
\begin{align}
 r^{-p}\Rea \mathcal C(z,h)&=(q-1)\operatorname{Var}_\mu(a)
       +\operatorname{Var}_\mu(b)
       +(p-1)(\mathbb E_\mu a)^2+(\mathbb E_\mu b)^2\ge 0,
       \\ 
 r^{-p}\Ima \mathcal C(z,h)
 &=(q-2)\operatorname{Cov}_\mu(a,b)
       +(p-2)(\mathbb E_\mu a)(\mathbb E_\mu b).
\end{align}
Use inequalities
\begin{align*}
|(q-2)\operatorname{Cov}_\mu(a,b)|
 &\le c_q\bigl((q-1)\operatorname{Var}_\mu(a)
+\operatorname{Var}_\mu(b)\bigr),\\
|(p-2)(\mathbb E_\mu a)(\mathbb E_\mu b)|
 &\le c_p\bigl((p-1)(\mathbb E_\mu a)^2 +(\mathbb E_\mu b)^2\bigr)
\end{align*}
Therefore,
\begin{equation}\label{eq:lem:key}
    |\Ima \mathcal C(z,h)|\le c_{p,q}\Rea \mathcal C(z,h).
\end{equation}
Next, we prove \eqref{eq:umd-two-point-sector}. Fix $u,v\in \mathbb C^N$. Assume $u \ne v$ or there is nothing to prove. Put $h=u-v$, $z(t)=v+th$, where $0\le t\le 1$. Then $z(0)=v$, $z(1)=u$ and $z'(t)=h$. At a parameter value where all the coordinates of $z(t)$ are nonzero, the chain rule implies
\begin{equation*}
    \frac{d}{dt} J_{p,q}(z(t)) = D J_{p,q}(z(t)) h.
\end{equation*}
For a moment, suppose one can integrate this derivative along the whole line segment $z(t)$ for $0\le t\le 1$. Then the fundamental theorem of calculus yields
\begin{equation*}
    J_{p,q}(u) - J_{p,q}(v) = \int_0^1 D J_{p,q}(z(t)) h\, dt.
\end{equation*}
With $h=u-v$, the above identity shows
\begin{align*}
    \sum_j (u_j-v_j) \left( \overline{J_{p,q}(u)_j-J_{p,q}(v)_j}    \right) &= \int_0^1  \sum_j h_j \overline{\left(D J_{p,q}(z(t)) h\right)_j}\,dt\\
    &=\int_0^1 \mathcal C(z(t),h)\,dt.
\end{align*}
Clearly, the real part of the left side equals $\int_0^1 \Rea \mathcal C(z(t),h)\,dt\ge 0$. Moreover, \eqref{eq:umd-two-point-sector} follows directly from \eqref{eq:lem:key}. We now treat the case when coordinates vanish. First discard coordinates that vanish identically along the segment. Each remaining $z_j(t)=v_j+th_j$ is affine and can vanish at most once. There are two cases. Suppose that $z_j(t_0)=0$ but $z(t_0)\ne 0$. Then $z_j(t)=(t-t_0)h_j$. Therefore, $|z_j(t)|^{q-2}z_j(t) = |h_j|^{q-2} h_j (t-t_0) |t-t_0|^{q-2}$ with derivative $(q-1)|h_j|^{q-2}h_j |t-t_0|^{q-2}$. Since $q>1$, the integral $\int_{t_0-\delta}^{t_0+\delta} |t-t_0|^{q-2}\,dt$ is finite and hence the derivative of $J_{p,q}(z(t))_j$ is integrable. If instead $z(t_0)=0$, then $v=-t_0 h$ and $z(t) = (t-t_0)h$. It follows that $J_{p,q}(z(t))= (t-t_0) |t-t_0|^{p-2} J_{p,q}(h)$, whose derivative $(p-1)|t-t_0|^{p-2} J_{p,q}(h)$ is again integrable near $t_0$ since $p>1$. The verification of \eqref{eq:umd-two-point-sector} is now complete.

Finally, we show sharpness of the constant $c_{p,q}$. Suppose that for some $c>0$, \eqref{eq:umd-two-point-sector} holds. Fix $z$ with nonzero coordinates and a direction $h$. Choose $v=z$ and $u=z+\varepsilon h$. One checks by differentiability that
\begin{equation*}
    J_{p,q}(z+\varepsilon h) - J_{p,q}(z) = \varepsilon D J_{p,q}(z) h + o(\varepsilon).
\end{equation*}
Thus 
\begin{equation*}
    \sum_j (u_j-v_j) \left( \overline{J_{p,q}(u)_j - J_{p,q}(v)_j} \right) = \varepsilon^2 \mathcal C(z,h) + o(\varepsilon^2).
\end{equation*}
Substituting this into \eqref{eq:umd-two-point-sector}, dividing both sides by $\varepsilon^2$ and then letting $\varepsilon \to 0$, one gets
\begin{equation*}
    |\Ima \mathcal C(z,h)| \le c \Rea \mathcal C(z,h) \iff c \ge \frac{|\Ima \mathcal C(z,h)|}{\Rea \mathcal C(z,h)}.
\end{equation*}
Take $a_j=1$, $b_j=\sqrt{p-1}$ for all $1\le j\le N$. One checks easily $\mathrm{Var}_\mu(a)=\mathrm{Var}_\mu(b)=\mathrm{Cov}_\mu(a,b)=0$ and $\mathbb E_\mu(a)=1$, $\mathbb E_\mu(b)=\sqrt{p-1}$. Immediately, $c\ge \frac{|p-2|}{2\sqrt{p-1}}=c_p$. For the $q$ obstruction, work first in two coordinates and take $z=(1,1)$, $a=(1,-1)$ and $b=\sqrt{q-1}(1,-1)$. The same test embeds into every $N\ge2$ by adjoining zero coordinates. Obviously, the weights are $\mu_1=\mu_2=1/2$ and $\mathbb E_\mu(a)=\mathbb E_\mu(b)=0$. A direct computation gives $\mathrm{Var}_\mu(a)=1$, $\mathrm{Var}_\mu(b)=q-1$, $\mathrm{Cov}_\mu(a,b)=\sqrt{q-1}$. Consequently, $c\ge \frac{|q-2|}{2\sqrt{q-1}}=c_q$. The proof of Lemma~\ref{lem:umd-mixed-sector} is now complete.
\end{proof}

The next lemma records the norm differentiation formula.

\begin{lemma}\label{lem:norm-derivative}
Let $E$ be uniformly smooth and let $1<p<\infty$. Recall the pointwise duality map $J_{p,E}$ with the
normalization
\[
 \|J_{p,E}(a)\|_{E^*}=\|a\|_E^{p-1},
 \qquad
 \langle a,J_{p,E}(a)\rangle=\|a\|_E^p,
 \qquad J_{p,E}(0)=0.
\]
Suppose that an $L^p(X;E)$-valued curve $v(t)$ satisfies
\[
 \frac{v(t)-v(t_0)}{t-t_0}\longrightarrow z
 \quad\hbox{in }L^p(X;E).
\]
Then
\begin{equation}\label{eq:norm-curve}
 \left.\frac d{dt}\|v(t)\|_{L^p(X;E)}^p\right|_{t=t_0}
 =
 p\Rea\int_X
 \langle z,J_{p,E}(v(t_0))\rangle\,d\mu.
\end{equation}
The same conclusion holds for the right derivative
when the hypothesis is assumed only as $t\downarrow t_0$.
\end{lemma}

\begin{proof}
This is the classical differentiation formula for a smooth Banach
norm; we include a proof in the present setting for convenience.
We use the fact that uniform smoothness of
$E$ implies norm continuity of the pointwise map
\[
 J_{p,E}:L^p(X;E)\longrightarrow L^{p'}(X;E^*),
 \qquad p'=\frac p{p-1};
\]
see \cite[Chapter~I, Sections~3--4, and Chapter~II, Section~4]{Cioranescu}. First, for $a,b\in E$, the norming properties of $J_{p,E}$ and convexity of $s\mapsto s^p/p$ give
\begin{equation*}
\Rea\langle b-a,J_{p,E}(a)\rangle
 \le \|a\|_E^{p-1}(\|b\|_E-\|a\|_E)\le \frac{\|b\|_E^p-\|a\|_E^p}{p}.
\end{equation*}
These inequalities also hold when $a=0$.
Interchanging $a$ and $b$ therefore yields
\[
 \Rea\langle b-a,J_{p,E}(a)\rangle
 \le \frac{\|b\|_E^p-\|a\|_E^p}{p}
 \le \Rea\langle b-a,J_{p,E}(b)\rangle.
\]
Apply this pointwise with $a=v(x)$ and $b=v(x)+h(x)$,
where $v,h\in L^p(X;E)$, and integrate. Subtracting
the first pairing and using H\"older's inequality, we obtain
\[
\begin{split}
 0
 &\le
 \frac{\|v+h\|_{L^p(X;E)}^p-\|v\|_{L^p(X;E)}^p}{p}
 -\Rea\int_X\langle h,J_{p,E}(v)\rangle\,d\mu\\
 &\le
 \|h\|_{L^p(X;E)}
 \|J_{p,E}(v+h)-J_{p,E}(v)\|_{L^{p'}(X;E^*)}.
\end{split}
\]
By continuity of the duality map, the last factor
tends to zero as $\|h\|_{L^p(X;E)}\to0$.

Now take $v=v(t_0)$ and $h=v(t)-v(t_0)$. The hypothesis
means that
\[
 h=(t-t_0)z+r(t),
 \qquad
 \frac{\|r(t)\|_{L^p(X;E)}}{|t-t_0|}\longrightarrow0.
\]
In particular, $\|h\|_{L^p(X;E)}=O(|t-t_0|)$.
The preceding remainder estimate takes the form
\[
 \frac{\|v(t)\|_{L^p(X;E)}^p
       -\|v(t_0)\|_{L^p(X;E)}^p}{p(t-t_0)}
 =
 \Rea\int_X
 \left\langle\frac{h}{t-t_0},
 J_{p,E}(v(t_0))\right\rangle\,d\mu+o(1).
\]
Since $h/(t-t_0)\to z$ in $L^p(X;E)$ and
$J_{p,E}(v(t_0))\in L^{p'}(X;E^*)$, H\"older's
inequality identifies the limit of the pairing.
This proves \eqref{eq:norm-curve}.
Every step remains valid with only $t\downarrow t_0$,
which proves the right-derivative assertion.
\end{proof}

We will also use the following classical Gaussian argument underlying the Marcinkiewicz--Zygmund theorem \cite{MarcinkiewiczZygmund} to supply Hilbert-valued extensions. We give its
mixed-norm form below for convenience.

\begin{lemma}\label{lem:gaussian-amplification}
Let $1<p,q<\infty$ and let $T$ be a bounded complex-linear
operator on $L^p(X)$. Suppose that
\[
 \|(Tf_j)_{j=1}^N\|_{L^p(X;\ell_N^q)}
 \le C\|(f_j)_{j=1}^N\|_{L^p(X;\ell_N^q)}
\]
for every $N\ge1$, with the same constant $C$.
Then, for every sigma-finite measure space $(Y,\nu)$ and every
separable complex Hilbert space $\mathcal{H}$, its action in the $X$ variable
extends to $L^p(X;L^q(Y;\mathcal{H}))$ and
\[
 \|TF\|_{L^p(X;L^q(Y;\mathcal{H}))}
 \le C\|F\|_{L^p(X;L^q(Y;\mathcal{H}))}.
\]
In particular, for every bounded complex-linear operator $T$ on
$L^r(X)$, $1<r<\infty$, its Hilbert-valued extension satisfies
\[
 \|T\|_{L^r(X;\mathcal{H})\to L^r(X;\mathcal{H})}
 =\|T\|_{L^r(X)\to L^r(X)},
\]
provided $\mathcal{H}\ne\{0\}$. All constants are independent of the
dimensions of the sequence and Hilbert variables.
\end{lemma}

\begin{proof}
First, the assumed estimate extends to every scalar inner
$L^q(Y)$ space. Indeed, for pairwise disjoint measurable sets
$A_j$ of finite measure, the norm of
$\sum_j f_j(x)\mathbf1_{A_j}(y)$ is the mixed sequence norm of
$(\nu(A_j)^{1/q}f_j(x))_j$. Applying the hypothesis to these
weighted coordinates and then using density proves this assertion.

For $\mathcal{H}=\mathbb C^m$, let $\gamma_1,\ldots,\gamma_m$ be independent
standard circular complex Gaussian variables on a probability
space $(\Omega_{\mathrm G},\mathbb P_{\mathrm G})$.
Thus $\gamma_k=(g_k+ih_k)/\sqrt2$, where the $g_k,h_k$ are
independent real $N(0,1)$ variables. The rotational invariance
of their joint distribution implies
\[
 \mathbb E\left|\sum_{k=1}^m\gamma_k z_k\right|^q
 =d_q\left(\sum_{k=1}^m|z_k|^2\right)^{q/2},
 \qquad d_q=\mathbb E|\gamma_1|^q.
\]
For $F=\sum_{k=1}^mF_ke_k$, where $(e_k)$ is an orthonormal
basis, put
\[
 \mathcal GF(x,y,\gamma)
 =\sum_{k=1}^m\gamma_kF_k(x,y).
\]
The moment identity, integrated in $y$, gives pointwise in $x$
\[
 \|\mathcal GF(x,\cdot,\cdot)\|_{L^q(Y\times\Omega_{\mathrm G},
                  \,d\nu\,d\mathbb P_{\mathrm G})}
 =d_q^{1/q}\|F(x,\cdot)\|_{L^q(Y;\mathcal{H})}.
\]
Since $T$ acts in $x$ alone, linearity gives
$T\mathcal GF=\mathcal G(TF)$. Apply the scalar inner-space
estimate to $Y\times\Omega_{\mathrm G}$. Taking the outer
$L^p(X)$ norm in the preceding identity, we obtain
\[
 \begin{split}
 d_q^{1/q}\|TF\|_{L^p(X;L^q(Y;\mathcal{H}))}
 &=\|T\mathcal GF\|_{L^p(X;L^q(Y\times\Omega_{\mathrm G}))}\\
 &\le C\|\mathcal GF\|_{L^p(X;L^q(Y\times\Omega_{\mathrm G}))}\\
 &=Cd_q^{1/q}\|F\|_{L^p(X;L^q(Y;\mathcal{H}))}.
 \end{split}
\]
The factor $d_q^{1/q}$ cancels. Orthogonal finite-rank
approximation proves the assertion for separable $\mathcal{H}$.

Finally, when $p=q=r$, Fubini verifies the hypothesis with
$C=\|T\|_{L^r\to L^r}$. Taking $Y$ to be a one-point space
proves the Hilbert-valued upper bound. Testing on $F(x)=f(x)e$
for a fixed unit vector $e\in \mathcal{H}$ proves the reverse inequality.
\end{proof}

We now prove complex contractivity.

\begin{proposition}\label{prop:umd-complex-contractivity}
Let $\mathcal{L}\ge0$ be selfadjoint on $L^2(X,\mu)$, and assume that $P_t=e^{-t\mathcal{L}}$ is a symmetric contraction semigroup on $L^r$ spaces, $1\le r\le\infty$. For $\Rea z\ge0$, define $P_z=e^{-z\mathcal{L}}$ first by the spectral theorem. For every separable complex Hilbert space $\mathcal{H}$, every $N\ge1$, and
$|\theta|\le\phi_{p,q}$,
\begin{equation}\label{eq:umd-complex-contractivity}
 \|P_{te^{i\theta}}F\|_{L^p(X;\ell^q_N(\mathcal{H}))}
 \le\|F\|_{L^p(X;\ell^q_N(\mathcal{H}))},\qquad t>0.
\end{equation}
The constants are independent of $N$ and $\mathcal{H}$. In particular the same assertion holds for $\ell^q(\mathcal{H})$.

If $0<\beta<1$ and
$V\in\operatorname{Dom}_{L^p(X;\ell^q(\mathcal{H}))}\mathcal{L}^{\beta}$, then
\begin{equation}\label{eq:umd-B-sector}
 |\Ima\langle \mathcal{L}^{\beta}V,J_{p,E}(V)\rangle|
 \le c_{p,q}\Rea\langle \mathcal{L}^{\beta}V,J_{p,E}(V)\rangle,
\end{equation}
where $E=\ell^q(\mathcal{H})$, the pairing in \eqref{eq:umd-B-sector} is integrated over $X$, and $J_{p,E}$ is its pointwise duality map, normalized
by $\langle V(x),J_{p,E}(V(x))\rangle_E=\|V(x)\|_E^p$.
\end{proposition}

\begin{proof}
We first take $\mathcal{H}=\mathbb C$. For $f=(f_1,\dots,f_N)$, we claim that
\begin{equation}\label{eq:prop-1}
\Rea\left(e^{i\theta}\sum_{j=1}^N \int_X
           (\mathcal{L}f_j)\overline{J_{p,q}(f)_j}\,d\mu\right)\ge0,
\end{equation}
whenever $f_j\in\operatorname{Dom}_{L^p}\mathcal{L}$. Equivalently, \eqref{eq:prop-1} is
\begin{equation*}
\Rea \sum_{j=1}^N \int_X
           \mathcal{L}\left(F_j(f)\right) G_j(f)\,d\mu \ge0,
\end{equation*}
with $F_j(f)=f_j$ and $G_j(f)=e^{i \theta} \overline{J_{p,q}(f)_j}$. Since $N$ is finite, one checks $f\in L^p(X;\ell_N^q)$. Moreover,
\[
 |J_{p,q}(f)_j|\le\|J_{p,q}(f)\|_{\ell_N^{q'}}
 =\|f\|_{\ell^q}^{p-1}.
\]
So, $\int |G_j(f)|^{p'}\,d\mu \le \|f\|_{L^p(X;\ell^q)}^p<\infty$. By \cite[Theorem~2.2]{HaaseForms}, it therefore suffices to prove \eqref{eq:prop-1} on the two-point probability space $X=\{0,1\}$ with measure $\mu(\{0\})=\mu(\{1\})=1/2$, replacing $\mathcal{L}$ by $I-E_\lambda$, where 
\[
E_\lambda = 
\begin{pmatrix}
0 & \lambda  \\
\overline{\lambda} & 0 
\end{pmatrix}, \qquad \forall|\lambda|=1.
\]
Assume $f(0)=u$, $f(1)=v$ and $u,v\in \mathbb C^N$. For each $1\le j\le N$, set $f_j = (u_j,v_j)^T$. Then proving \eqref{eq:prop-1} is reduced to showing
\begin{equation*}
    \Rea \left( e^{i \theta} \mathcal Q_\lambda(u,v)\right) \ge 0,\qquad \forall u,v\in \mathbb C^N,\qquad \forall |\lambda|=1,
\end{equation*}
where
\begin{equation*}
\mathcal Q_\lambda(u,v) = \frac{1}{2} \sum_j (u_j-\lambda v_j) \overline{J_{p,q}(u)_j} + \frac{1}{2} \sum_j (v_j-\overline{\lambda}u_j) \overline{J_{p,q}(v)_j}.
\end{equation*}

We first show that it holds for $\lambda=1$. In this case,
\begin{equation*}
\mathcal Q_1(u,v) = \frac{1}{2} \sum_j (u_j-v_j) \left( \overline{J_{p,q}(u)_j - J_{p,q}(v)_j}  \right).
\end{equation*}
From Lemma~\ref{lem:umd-mixed-sector}, one has $\Rea \mathcal{Q}_1 \ge 0$ and $|\Ima \mathcal Q_1|\le c_{p,q} \Rea \mathcal{Q}_1$. Consequently,
\begin{align*}
\Rea \left( e^{i \theta} \mathcal Q_1\right) &= \cos{\theta} \,\Rea \mathcal{Q}_1 - \sin{\theta}\, \Ima \mathcal Q_1\\ 
&\ge \left( \cos{\theta} - c_{p,q}|\sin{\theta}|\right) \Rea \mathcal{Q}_1 \ge 0
\end{align*}
since $c_{p,q}=\cot{\phi_{p,q}}$ and $|\theta|\le \phi_{p,q}$.

Next, for general $|\lambda|=1$, the definition of $J_{p,q}$ gives $J_{p,q}(\lambda \cdot) = \lambda J_{p,q}(\cdot)$. Write $w = \lambda v$. One obtains 
\begin{align*}
\Rea \left(e^{i\theta}\mathcal Q_\lambda(u,v) \right) &= \Rea \left(e^{i\theta} \frac{1}{2} \sum_j (u_j-w_j) \left( \overline{J_{p,q}(u)_j-J_{p,q}(w)_j} \right) \right) \\
&= \Rea \left(e^{i\theta} \mathcal Q_1(u,w) \right)\ge 0.
\end{align*}
This confirms the claim i.e., \eqref{eq:prop-1}.

For finite $N$, the scalar holomorphic semigroup, supplied by \cite[Corollary~6.1]{HaaseForms}, acts coordinatewise on $L^p(X;\mathbb C^N)$. On each ray strictly inside the scalar contractivity sector $|\theta|<\phi_{p,q}\le \phi_p$, its trajectories are differentiable. Apply Lemma~\ref{lem:norm-derivative} with $E=\ell_N^q$ to the curve $v(t)=P_{te^{i\theta}}f$. For $t>0$,
its derivative in $L^p(X;\ell_N^q)$ is $v'(t)=-e^{i\theta}\mathcal{L}P_{te^{i\theta}}f$. The pointwise duality map $J_{p,E}$ is precisely $J_{p,q}$ from Lemma~\ref{lem:umd-mixed-sector}.
Consequently, \eqref{eq:norm-curve} becomes
\begin{equation*}
\frac d{dt}\|P_{te^{i\theta}}f\|_{L^p(\ell^q_N)}^p
 =-p\Rea\left(e^{i\theta}\sum_j\int_X \mathcal{L}P_{te^{i\theta}}f_j\, \overline{J_{p,q}(P_{te^{i\theta}}f)_j}\,d\mu\right) \le 0.
\end{equation*}
Strong continuity implies \eqref{eq:umd-complex-contractivity}. 

Apply Lemma~\ref{lem:gaussian-amplification} to
$T=P_{te^{i\theta}}$ with $C=1$, using the finite-coordinate
estimate just proved. Taking the inner measure space to be
counting measure on $\{1,\ldots,N\}$, or on $\mathbb N$,
proves the assertion for $\ell_N^q(\mathcal{H})$ and $\ell^q(\mathcal{H})$.

For the boundary rays, first take finite-coordinate, finite
Hilbert-rank functions whose scalar entries lie in $L^2\cap L^p$.
As interior angles approach either boundary, the spectral theorem
ensures strong $L^2$ convergence to the boundary operator.
An almost everywhere convergent subsequence and Fatou's lemma
preserve the contraction bound. Density and finite-rank
approximation complete the passage to $L^p(X;\ell^q(\mathcal{H}))$.
When $\phi_{p,q}=\pi/2$, necessarily $p=q=2$, and the
Hilbert spectral theorem supplies the imaginary boundary directly.

Finally, Bochner subordination \cite{Stein} makes $\mathcal{L}^\beta$ a symmetric
contraction generator. Apply the established complex
contractivity to this generator. For $E=\ell^q(\mathcal{H})$, $V\in\operatorname{Dom}_{L^p(X;E)}\mathcal{L}^{\beta}$ and
$|\theta|<\phi_{p,q}$, set $v(t)=e^{-te^{i\theta}\mathcal{L}^{\beta}}V$ for $t\ge 0$. The generator identity reads
\[
 \frac{v(t)-V}{t}\longrightarrow-e^{i\theta}\mathcal{L}^{\beta}V
 \quad\hbox{in }L^p(X;E)
\]
as $t \to 0$. By complex contractivity, $\|v(t)\|_{L^p(X;E)}\le\|V\|_{L^p(X;E)}$; the right-derivative formula in
Lemma~\ref{lem:norm-derivative} yields
\[
 0\ge
 \lim_{t\downarrow0}
 \frac{\|v(t)\|_{L^p(X;E)}^p-\|V\|_{L^p(X;E)}^p}{pt}
 =
 -\Rea\left(e^{i\theta}
       \langle \mathcal{L}^{\beta}V,J_{p,E}(V)\rangle\right).
\]
Equivalently,
\begin{equation*}
 \Rea\left(e^{i\theta}
       \langle \mathcal{L}^{\beta}V,J_{p,E}(V)\rangle\right)\ge0,
 \qquad |\theta|<\phi_{p,q}.
\end{equation*}
Expanding this inequality and letting $|\theta| \to \phi_{p,q}$ completes the proof of \eqref{eq:umd-B-sector} and hence Proposition~\ref{prop:umd-complex-contractivity}.
    
\end{proof}

\subsection{A contact inequality}\label{subsec:contact}

Throughout this subsection, we assume
\begin{equation}
1<p<\min(q,2).
\end{equation}
Fix a finite number $n$ of coordinates and a compact interval
$I\subset(0,\infty)$. Put
\[
 \mathcal{H}=L^2(I,ds/s),\qquad E=\ell^q_n(\mathcal{H}),\qquad
 (Kv)_j(s)=s^{-1/2}v_j(s).
\]
The operators $\mathcal{L}^{1/2}$ and $K$ commute, and $K$ is bounded and boundedly invertible on the current scale interval.

Write $|\cdot|_{\mathcal{H}}$ for the Hilbert norm. For $E=\ell_n^q(\mathcal{H})$ we also write $J_{p,q}=J_{p,E}$; explicitly,
\[
 J_{p,q}(v)_j=\|v\|_E^{p-q}\|v_j\|_{\mathcal{H}}^{q-2}v_j,
 \qquad E^*=\ell_n^{q'}(\mathcal{H}),
\]
with value zero in each zero coordinate.

\begin{proposition}
\label{prop:common-contact}
Let $f\in L^p(X;E)\cap L^q(X;E)\cap L^2(X;E)$ and $h>0$.
There are $u\in\bigcap_{s\in\{p,q,2\}}\operatorname{Dom}_{L^s(X;E)}\mathcal{L}^{1/2}$ and $g\in\bigcap_{s\in\{p,q,2\}}L^s(X;E)$ with $f=(\mathcal{L}^{1/2}+K)u+g$, $\|g(x)\|_E\le h$,
and $g=hu/\|u\|_E$ on $\Omega=\{u\ne0\}$, satisfying
\eqref{um:obstacle-budget}. Put $V=K^{1/p'}u$ and $Z=\Omega^c$.
Then
\begin{equation}\label{eq:energy}
0\le\mathcal E(V):=\Rea\langle \mathcal{L}^{1/2}V,J_{p,q}(V)\rangle
 \le \frac{(p-1)^{p-1}}{p^p}\|f\|_{L^p(E)}^p,
\end{equation}
and
\begin{equation}\label{eq:contact}
\begin{aligned}
 \int_0^\infty
 \|\one_Ze^{-te^{i\theta}\mathcal{L}}V\|_{L^p(E)}^p\frac{dt}{t^{3/2}}
 &\le 2p\sqrt{\pi}\,
       \Rea\left(e^{i\theta/2}\langle \mathcal{L}^{1/2}V,J_{p,q}(V)\rangle\right)\\
 &\le 2p\sqrt{\pi}(1+c_{p,q})\mathcal E(V)\\
 &\le C_{p,q}\|f\|_{L^p(E)}^p,
 \qquad |\theta|\le\phi_{p,q}.
\end{aligned}
\end{equation}
The constants are independent of $n$ and $I$.
\end{proposition}

\begin{proof}
We apply Theorem~\ref{lem:radial} with $\mathcal{A}=\mathcal{L}^{1/2}+K$. Indeed, write $\kappa:=\inf_{s\in I}s^{-1/2}$. Since $I\subset (0,\infty)$ is compact, $\kappa>0$, and
\[
 \|e^{-t(\mathcal{L}^{1/2}+K)}\|_{L^{r}(X;E)\to L^r(X;E)}
 \le e^{-\kappa t}
\]
for every $1\le r\le \infty$. For each $t$, let $S_t$ be the positive Dunford--Schwartz modulus of the operator $e^{-t\mathcal{L}^{1/2}}$ such that $S_t \mathbf{1}\le \mathbf{1}$; its existence and domination are classical, see \cite[Lemma~3.4]{Taggart}. Therefore,
\begin{equation*}
\|e^{-t(\mathcal{L}^{1/2}+K)}u(x)\|_E\le S_t(\|u\|_E)(x).
\end{equation*}
Next, let $\Psi$ be convex and increasing, with $\Psi(0)=0$. Monotonicity and Jensen's inequality \cite[Theorem~3.4 (b)]{HaaseConvexity} give
\begin{equation}\label{eq:Jensen}
\Psi\left( \|e^{-t(\mathcal{L}^{1/2}+K)}u(x)\|_E \right) \le \Psi\left( S_t(\|u\|_E)(x) \right) \le S_t \Psi(\|u\|_E)(x).
\end{equation}
Indeed, $\Psi$ is the supremum of supporting affine functions $\ell(s)=c_1s+c_2$ with $c_1\ge0$ and $c_2\le0$. For each such function, $\Psi(s)\ge\ell(s)$ for all $s\ge0$. Positivity and $S_t \mathbf{1}\le \mathbf{1}$ imply
\begin{equation*}
    S_t \Psi(\|u\|_E) \ge c_1 S_t(\|u\|_E) + c_2S_t \mathbf{1} \ge c_1 S_t(\|u\|_E) + c_2 = \ell(S_t(\|u\|_E)).
\end{equation*}
Taking the supremum over all supporting functions shows \eqref{eq:Jensen}. Integration and $L^1$ contractivity verify the remaining hypothesis of Theorem~\ref{lem:radial}. Moreover, $E$ and $L^r(X;E)$ are uniformly smooth and reflexive for $1<r<\infty$. Run the fixed-point construction in the proof of
Theorem~\ref{lem:radial} on
$L^p(X;E)\cap L^q(X;E)\cap L^2(X;E)$ with the sum norm.
The contraction factor is the same in all three spaces, and
$\|u_\varepsilon\|_{L^s(E)}\le\kappa^{-1}\|f\|_{L^s(E)}$
for $s=p,q,2$.
The pointwise cap and the $L^p$ estimate ensure
$\|g_\varepsilon\|_{L^s(E)}^s
\le h^{s-p}\|f\|_{L^p(E)}^p$ for $s=q,2$.
The same strong $L^p$ limit and a common weak subsequence
therefore produce $u$ in all three generator domains.
Thus we can decompose
\begin{equation}\label{um:obstacle}
 f=(\mathcal{L}^{1/2}+K)u+g,\qquad \|g(x)\|_E\le h,\qquad
 g=h\frac u{\|u\|_E}\ \text{on }\Omega=\{u\ne0\},
\end{equation}
with
\begin{equation}\label{um:obstacle-budget}
 \|g\|_{L^p(E)}\le\|f\|_{L^p(E)},\quad
 \mu(\Omega)\le h^{-p}\|f\|_{L^p(E)}^p,\quad
 \|g\|_{L^q(E)}^q\le h^{q-p}\|f\|_{L^p(E)}^p.
\end{equation}
Note that since $K$ is bounded, $u\in\operatorname{Dom}_{L^s}(\mathcal{L}^{1/2})$ for $s=2,p,q$. Define
\begin{equation*}
V=K^{1/p'}u,\qquad Z=\Omega^c.
\end{equation*}
Since $K$ is invertible, $V$ vanishes on $Z$.

Recall the duality map 
\begin{equation*}
J_{p,q}(v)_j=\|v\|_E^{p-q}|v_j|_{\mathcal{H}}^{q-2}v_j,\qquad J_{p,q}(v)\in E^*=\ell_n^{q'}(\mathcal{H}),
\end{equation*}
with the usual zero convention and set
\begin{equation*}
Q(V)=\langle KV,J_{p,q}(V)\rangle_E
 =\|V\|_E^{p-q}\sum_j |V_j|_{\mathcal{H}}^{q-2}|K^{1/2}V_j|_{\mathcal{H}}^2.
\end{equation*}
Note that $Q$ is a nonnegative real-valued function of $x$.

We claim that
\begin{equation}\label{um:killing-holder}
       \|K^{1/p'}J_{p,q}(V)\|_{E^*}\le Q(V)^{1/p'}.
\end{equation}
If $V=0$ then both sides vanish. Suppose $V\ne0$, put $\nu_j= |V_j|_{\mathcal{H}}^q/\|V\|_E^q \ge 0$. Then $\sum_j\nu_j=1$, so these numbers are probability weights. Also, for each nonzero coordinate (zero coordinates can be omitted), define
\begin{equation*}
r_j=\frac{|K^{1/2}V_j|_{\mathcal{H}}^2}{|V_j|_{\mathcal{H}}^2}.
\end{equation*}
It is plain that $Q(V)=\|V\|_E^{p-q} \sum_j |V_j|_{\mathcal{H}}^q r_j = \|V\|_E^p \sum_j \nu_jr_j$.

Since $1/p'<1/2$, H\"older's inequality gives 
\begin{align*}
|K^{1/p'}V_j|_{\mathcal{H}}^2 &= \int_I s^{-1/p'} |V_j(s)|^2\, \frac{ds}{s}\\
&\le \left( \int_I s^{-1/2} |V_j(s)|^2\,\frac{ds}{s} \right)^{2/p'} \left( \int_I |V_j(s)|^2\,\frac{ds}{s} \right)^{1-2/p'}\\
&= |K^{1/2} V_j|_{\mathcal{H}}^{4/p'} |V_j|_{\mathcal{H}}^{2-4/p'}.
\end{align*}
That is, $|K^{1/p'}V_j|_{\mathcal{H}}\le |V_j|_{\mathcal{H}} r_j^{1/p'}$. Consequently,
\begin{align*}
\|K^{1/p'}J_{p,q}(V)\|_{E^*}
 &= \|V\|_E^{p-q}
       \left(\sum_j |V_j|_{\mathcal{H}}^{(q-2)q'}
       |K^{1/p'}V_j|_{\mathcal{H}}^{q'}\right)^{1/q'}\\
 &\le \|V\|_E^{p-q}
       \left(\sum_j |V_j|_{\mathcal{H}}^q r_j^{q'/p'}\right)^{1/q'}\\
 &= \|V\|_E^{p-1}
       \left(\sum_j\nu_j r_j^{q'/p'}\right)^{1/q'}\\
 &\le \|V\|_E^{p-1}
       \left(\sum_j\nu_jr_j\right)^{1/p'}\\
 &= Q(V)^{1/p'},
\end{align*}
where the last inequality follows by H\"older's inequality with exponent $p'/q'\ge 1$ and $\sum_j \nu_j=1$. This proves the claim \eqref{um:killing-holder}.

Multiplying \eqref{um:obstacle} by $K^{1/p'}$, we find $\mathcal{L}^{1/2}V+KV=K^{1/p'}(f-g)$. Pair it with $J_{p,q}(V)$, integrate over $X$, and take real parts:
\begin{equation*}
\mathcal E(V) + \int_X Q(V) \,d\mu = \Rea \left \langle K^{1/p'}f,\, J_{p,q}(V)    \right \rangle - \Rea \left \langle K^{1/p'}g,\, J_{p,q}(V)    \right \rangle
\end{equation*}
On $\Omega$, it is clear that $K^{1/p'}g=hV/\|u\|_E$, so
\[
 \left\langle K^{1/p'}g,\,J_{p,q}(V)\right\rangle_E
 =h\|V\|_E^p/\|u\|_E\ge0.
\] On $Z$, since $J_{p,q}(V)=0$, the pairing vanishes as well. By \eqref{um:killing-holder}, one deduces
\begin{equation}\label{eq:mathcalE}
\mathcal E(V) + \int_X Q(V) \,d\mu \le \int_X \|f(x)\|_E Q(V)(x)^{1/p'}\,d\mu.
\end{equation}
Observe that since $e^{-t\mathcal{L}^{1/2}}$ is contractive on $L^p(X;E)$, its right norm derivative near zero is non-positive, i.e., $\frac{d}{dt} \|e^{-t\mathcal{L}^{1/2}}V\|_{L^p(X;E)}^p \bigg|_{t=0^+}\le 0$. It follows by Lemma~\ref{lem:norm-derivative} that 
\begin{align*}
\frac{d}{dt} \|e^{-t\mathcal{L}^{1/2}}V\|_{L^p(X;E)}^p \bigg|_{t=0^+} = -p \Rea \left\langle \mathcal{L}^{1/2}V,\, J_{p,q}(V) \right\rangle = -p \mathcal E(V) \le 0.
\end{align*}
Hence $\mathcal E(V)\ge 0$.

Moreover, from \eqref{eq:mathcalE},
\begin{equation*}
\mathcal E(V) \le \int_X \left(\|f(x)\|_E Q(V)(x)^{1/p'} - Q(V)(x)\right) \,d\mu.    
\end{equation*}
Optimizing the function $\phi(s)= C s^{1/p'}-s$ with $C,s\ge 0$ yields $\phi(s)\le p^{-p} (p-1)^{p-1} C^p$. Hence
\begin{equation*}
\mathcal E(V)\le \frac{(p-1)^{p-1}}{p^p}\|f\|_{L^p(X;E)}^p.
\end{equation*}
This proves \eqref{eq:energy}.

Finally, we prove \eqref{eq:contact}. Let $|\theta|\le\phi_{p,q}$ and write $W=e^{-te^{i\theta}\mathcal{L}}V$. Proposition~\ref{prop:umd-complex-contractivity} gives $\|W\|_{L^p(X;E)}^p \le \|V\|_{L^p(X;E)}^p$. Clearly, for $x\in \Omega$,
\begin{align*}
    \Rea \langle W(x)-V(x),\, J_{p,q}(V(x))\rangle_E &= \Rea \langle W(x),\, J_{p,q}(V(x))\rangle_E - \|V(x)\|_E^p\\
    &\le \|W(x)\|_E \|J_{p,q}(V(x))\|_{E^*} - \|V(x)\|_E^p\\
    &= \|V(x)\|_E^{p-1} \left(  \|W(x)\|_E - \|V(x)\|_E \right).
\end{align*}
Convexity thus implies
\begin{equation*}
    \|W(x)\|_E^p \ge \|V(x)\|_E^p +p \Rea \left\langle W(x)-V(x), \, J_{p,q}(V(x)) \right\rangle_E.
\end{equation*}
On $Z$, both $V$ and $J_{p,q}(V)$ vanish. Integrating over $\Omega$ and adding the $Z$ contribution, we obtain
\begin{align*}
    \|W\|_{L^p(X;E)}^p &\ge \|V\|_{L^p(X;E)}^p + p \Rea \left\langle W-V, \, J_{p,q}(V) \right\rangle + \| \mathbf{1}_Z W \|_{L^p(X;E)}^p\\
    &= \|V\|_{L^p(X;E)}^p - p D(t) + \| \mathbf{1}_Z W \|_{L^p(X;E)}^p,
\end{align*}
where
\begin{equation*}
D(t):=\|V\|_{L^p(X;E)}^p-\Rea\langle W,J_{p,q}(V)\rangle.
\end{equation*}
Combining this with complex-contractivity of $W$, one gets
\begin{equation*}
    pD(t) \ge \| \mathbf{1}_Z W \|_{L^p(X;E)}^p \iff D(t) \ge \frac{1}{p} \| \mathbf{1}_Z e^{-te^{i\theta}\mathcal{L}}V\|_{L^p(X;E)}^p \ge 0.
\end{equation*}
Using the scalar identity
\begin{equation*}
e^{i\theta/2} \lambda^{1/2} = \frac{1}{2\sqrt{\pi}} \int_0^\infty \left(1-e^{-t e^{i\theta}\lambda} \right) \, \frac{dt}{t^{3/2}},
\end{equation*}
the spectral theorem gives the following integral
\begin{equation*}
e^{i\theta/2} \mathcal{L}^{1/2}V = \frac{1}{2\sqrt{\pi}} \int_0^\infty \left(V-e^{-t e^{i\theta}\mathcal{L}}V \right) \, \frac{dt}{t^{3/2}}.
\end{equation*}
Pairing with $J_{p,q}(V)$ and taking real parts, one obtains
\begin{equation*}
    \Rea \left( e^{i\theta/2} \langle \mathcal{L}^{1/2}V,\, J_{p,q}(V) \rangle \right) = \frac{1}{2\sqrt{\pi}} \int_0^\infty D(t)\, \frac{dt}{t^{3/2}}.
\end{equation*}
It follows that
\begin{align*}
    \frac{1}{2p\sqrt{\pi}} \int_0^\infty \|\mathbf1_Ze^{-te^{i\theta}\mathcal{L}}V\|_{L^p(X;E)}^p \, \frac{dt}{t^{3/2}} &\le \Rea \left( e^{i\theta/2} \langle \mathcal{L}^{1/2}V,\, J_{p,q}(V) \rangle \right)\\
    &\le (1+c_{p,q}) \mathcal E(V),
\end{align*}
where the last inequality follows by Proposition~\ref{prop:umd-complex-contractivity}.

We briefly justify the fractional identity on
$L^2(X;\ell_n^2(\mathcal{H}))$.
The simultaneous $L^2(X;E)$ construction and
finite-coordinate norm equivalence give
$V\in\operatorname{Dom}_{L^2(X;\ell_n^2(\mathcal{H}))}(\mathcal{L}^{1/2})$.
Since $J_{p,q}(V)$ vanishes outside $\Omega$,
H\"older's inequality and $1<p<2$ give
\[
 \|J_{p,q}(V)\|_{L^2(X;\ell_n^2(\mathcal{H}))}^2
 \le C_{n,q}^2
 \left(\int_\Omega\|V(x)\|_E^2\,d\mu\right)^{p-1}
 \mu(\Omega)^{2-p}<\infty.
\]
We may therefore pair the truncated Hilbert fractional
formula with $J_{p,q}(V)$. Strong convergence of the
truncated vectors, followed by monotone convergence
for the nonnegative defect $D(t)$, yields the asserted
identity and hence \eqref{eq:contact}.
The norm-equivalence constants establish membership
only and do not enter the contact estimate.

\end{proof}

\begin{remark}
The common obstacle is important for the independence of the
estimates from the number of coordinates. A coordinatewise
obstacle construction followed by a union bound would
produce the exceptional-set budget
\[
 h^{-p}\sum_{j=1}^n\|f_j\|_{L^p(X;\mathcal{H})}^p.
\]
In the range $p<q$, this quantity cannot be controlled
by $h^{-p}\|f\|_{L^p(X;\ell_n^q(\mathcal{H}))}^p$ with a constant
independent of $n$. Indeed, taking
$f=(f_0,\ldots,f_0)$ with $f_0\ne0$, one has
\[
 \sum_{j=1}^n\|f_j\|_{L^p(X;\mathcal{H})}^p
 =n\|f_0\|_{L^p(X;\mathcal{H})}^p
 =n^{1-p/q}\|f\|_{L^p(X;\ell_n^q(\mathcal{H}))}^p.
\]
Thus this direct coordinatewise estimate incurs a factor
$n^{1-p/q}$. Imposing the obstacle on the full fibre norm
gives the required mixed-norm budget directly. Moreover,
all coordinates and scales share the same exceptional set,
so the contact estimate can be applied to the entire fibre
before taking its norm.
\end{remark}

\subsection{Mellin reconstruction}\label{subsec:mellin}
The preceding proposition controls a complex semigroup trajectory
outside the obstacle support. We next express the spectral
multipliers used below as superpositions of these trajectories.

\begin{lemma}\label{lem:mixed-mellin}
Let $1<p<\min(q,2)$, $a=1/(2p)$, let $I\subset(0,\infty)$ be compact, and let $(s,\lambda)\mapsto H_s(\lambda)$ be jointly measurable, with $H_s:(0,\infty)\to\C$ continuous for each $s\in I$. Put $h_s(r)=H_s(e^r)$. Assume $h_s\in L^1(\R)$ for each $s\in I$. Suppose
\begin{equation}\label{um:mellin-size}
 \sup_{s\in I}
 \left\|(1+|\xi|)^{a+1/2}e^{\sigma_{p,q}|\xi|} \widehat h_s(\xi) \right\|_{L^p(d\xi)}\le M_0.
\end{equation}
Then there are two scalar functions $k_s^{\pm \phi_{p,q}}$ such that
\begin{align}
 \lambda^a H_s(\lambda)
 &=\sum_{\theta=\pm\phi_{p,q}}
     \int_0^\infty k_s^\theta(t)(e^{-te^{i\theta}\lambda}-1)\,dt,
\label{um:mellin-rep}\\ 
 \sup_{s\in I} \left\|t^{1+a}k_s^\theta(t) \right\|_{L^{p'}(dt/t)}
 &\le C_{p,q}M_0                                 \label{um:mellin-density}
\end{align}
for $\theta=\pm\phi_{p,q}$.
\end{lemma}

\begin{proof}
Choose $\theta=-\phi_{p,q}$ if $\xi>0$ and
$\theta=\phi_{p,q}$ if $\xi<0$. On the selected half-line,
$\theta\xi=-\phi_{p,q}|\xi|$. Let $\Gamma$ denote the Euler gamma function. Set
\[
 A_s^\theta(\xi)
 =\mathbf1_{\{\theta\xi<0\}}
   \frac{e^{\theta\xi}\widehat h_s(\xi)}
        {\Gamma(-a-i\xi)},
\]
and define
\[
 k_s^\theta(t)
 =\frac{e^{-ia\theta}}{2\pi}t^{-1-a}
   \widehat{A_s^\theta}(\log t).
\]
Here the Fourier transform of $A_s^\theta\in L^p(\mathbb R)$
is understood in the Hausdorff--Young sense.
We check that $\|A_s^\theta\|_p\le C_{p,q}M_0$.
Indeed, Stirling's formula shows
\[
 |\Gamma(-a-i\xi)|^{-1}
 \le C_a(1+|\xi|)^{a+1/2}e^{\pi|\xi|/2}.
\]
Thus, on the selected half-line,
\begin{align*}
 |A_s^\theta(\xi)|
 &\le C_a(1+|\xi|)^{a+1/2}
       e^{(\pi/2-\phi_{p,q})|\xi|}
       |\widehat h_s(\xi)|\\
 &=C_a(1+|\xi|)^{a+1/2}
       e^{\sigma_{p,q}|\xi|}
       |\widehat h_s(\xi)|.
\end{align*}
Taking the $L^p(\mathbb R)$ norm and using
\eqref{um:mellin-size} proves the claimed bound.
Hausdorff--Young and the change of variables $r=\log t$
therefore yield
\[
 \|t^{1+a}k_s^\theta(t)\|_{L^{p'}(dt/t)}
 =\frac1{2\pi}
   \|\widehat{A_s^\theta}\|_{L^{p'}(\mathbb R)}
 \le C_p\|A_s^\theta\|_p
 \le C_{p,q}M_0.
\]
This proves \eqref{um:mellin-density}.

Next, we verify \eqref{um:mellin-rep}. Integration by
parts and the gamma integral give
\[
 \int_0^\infty
 t^{-1-a-i\xi}(e^{-te^{i\theta}\lambda}-1)\,dt
 =
 \Gamma(-a-i\xi)e^{ia\theta-\theta\xi}
 \lambda^{a+i\xi},
 \qquad \lambda>0.
\]
The integral is absolutely convergent because $0<a<1$
and $|\theta|<\pi/2$.

The assumption \eqref{um:mellin-size} and H\"older's
inequality imply $\widehat h_s\in L^1(\mathbb R)$.
Since $h_s\in L^1(\mathbb R)$ is continuous, Fourier
inversion gives
\[
 \lambda^aH_s(\lambda)
 =\frac{\lambda^a}{2\pi}
   \int_{\mathbb R}
   \widehat h_s(\xi)e^{i\xi\log\lambda}\,d\xi.
\]

We justify the reconstruction directly by Fourier
duality in logarithmic time. For fixed $\lambda>0$,
the function $\phi(r) := e^{-ar}(e^{-\lambda e^{r+i\theta}}-1)$ belongs to the Schwartz space $\mathcal S(\mathbb R)$, since $0<a<1$ and $|\theta|<\pi/2$. Its Fourier transform is the gamma integral displayed above. Consequently, with the change of variables $t=e^r$,
\begin{align*}
 \int_0^\infty
 k_s^\theta(t)(e^{-te^{i\theta}\lambda}-1)\,dt &=  \frac{e^{-ia\theta}}{2\pi} \int_{\mathbb R} \widehat{A_s^\theta}(r) \phi(r)\,dr \\
 &=\frac{e^{-ia\theta}}{2\pi}
 \int_{\mathbb R} A_s^\theta(\xi)
 \left[
 \int_0^\infty
 t^{-1-a-i\xi}(e^{-te^{i\theta}\lambda}-1)\,dt
 \right]d\xi\\
 &=
 \frac{\lambda^a}{2\pi}
 \int_{\{\theta\xi<0\}}
 \widehat h_s(\xi)e^{i\xi\log\lambda}\,d\xi.
\end{align*}
Summing over $\theta=\pm\phi_{p,q}$ and applying
Fourier inversion proves \eqref{um:mellin-rep}.

For completeness, the time integrals satisfy the
uniform absolute bound
\begin{align*}
 \int_0^\infty
 |k_s^\theta(t)|\min(t\lambda,2)\,dt &\le
 \|t^{1+a}k_s^\theta\|_{L^{p'}(dt/t)}
 \left(
 \int_0^\infty
 t^{-ap}\min(t\lambda,2)^p\,\frac{dt}{t}
 \right)^{1/p}\\
 &\quad\le C_{p,q}M_0\lambda^a.
\end{align*}
The last inequality follows by substituting
$u=t\lambda$ and using $0<a<1$.
Since
$|e^{-te^{i\theta}\lambda}-1|\le\min(t\lambda,2)$,
this also proves absolute convergence of the
reconstruction integrals. 
\end{proof}

\begin{remark}
The exponential weight in \eqref{um:mellin-size} records the
angle of complex contractivity. In the definition of
$A_s^\theta$, Stirling's formula contributes the exponential
factor $e^{\pi|\xi|/2}$, while the selected half-line gives
$e^{\theta\xi}=e^{-\phi_{p,q}|\xi|}$. Their product is exactly
\[
 e^{\pi|\xi|/2}e^{-\phi_{p,q}|\xi|}
 =e^{\sigma_{p,q}|\xi|}.
\]
This explains the correspondence between the semigroup rays
and the Mellin weight. If one instead fixes interior rays of angles
$\pm(\phi_{p,q}-\varepsilon)$, where
$0<\varepsilon<\phi_{p,q}$, the same calculation produces
$e^{(\sigma_{p,q}+\varepsilon)|\xi|}$. The boundary rays
therefore match the stated hypothesis without an angular
loss in the present reconstruction.
\end{remark}

\begin{lemma}\label{lem:diagonal-time}
Let $1<q<\infty$ and $1<p\le\min(q,2)$ and
\[
 E=\ell^q_n(L^2(I,ds/s)).
\]
If $(s,t)\mapsto b_s(t)$ is jointly measurable and
$$\sup_s\|b_s\|_{L^{p'}(dt/t)}\le M_1$$ 
for some $M_1>0$, then for any $x(t)=\left( x_{j,s}(t) \right)_{j,s} \in L^p\left((0,\infty), dt/t; E\right)$,
\begin{equation}
 \left\|\left(\int_0^\infty b_s(t)x_{j,s}(t)\frac{dt}{t}
                                      \right)_{j,s}\right\|_E
 \le M_1\left(\int_0^\infty
               \|x(t)\|_E^p\frac{dt}{t}\right)^{1/p}.
\end{equation}
\end{lemma}

\begin{proof}
For clarity, note that $\|x(t)\|_E = \left( \sum_{j=1}^n \left( \int_I |x_{j,s}(t)|^2\, \frac{ds}{s} \right)^{q/2} \right)^{1/q}$ and $\int_0^\infty \|x(t)\|_E^p \,\frac{dt}{t}<\infty$. Define $y=(y_j(s))_{j,s}$, where
\begin{equation*}
    y_j(s) = \int_0^\infty b_s(t) x_{j,s}(t)\, \frac{dt}{t}.
\end{equation*}
Thus it is enough to show that $\|y\|_E \le M_1 \|x\|_{L^p((0,\infty),dt/t;E)}$. 

Put $d\nu(t)=dt/t$ and $d\rho(s)=ds/s$. By density, it suffices to consider simple functions supported on compact time intervals, say of the form $x(t) = \sum_{\ell=1}^N v_\ell \mathbf{1}_{A_\ell}(t)$, where $v_\ell \in E$ and each $A_\ell$ has finite $\nu$-measure. 

Let $w=(w_j(s))_j\in E^*=\ell^{q'}_n(L^2(I,d\rho))$ with $\|w\|_{E^*}\le1$. By duality, one only needs to prove $|\langle y, w\rangle|\le M_1 \|x\|_{L^p((0,\infty),d\nu;E)}$ for all such $w$. By Fubini's theorem
\begin{align*}
\langle y, w\rangle &= \sum_{j=1}^n \int_I y_j(s) \overline{w_j(s)} \, d\rho(s)\\
&= \sum_{j=1}^n \int_I \left( \int_0^\infty b_s(t) x_{j,s}(t) \, d\nu(t)  \right)  \overline{w_j(s)} \, d\rho(s)\\
&= \int_0^\infty \sum_{j=1}^n \int_I b_s(t) x_{j,s}(t) \overline{w_j(s)} \,d\rho(s) \,d\nu(t).
\end{align*}
For each $t>0$, define $F(t)= \left(F_j(t,s)\right)_j = \left( \overline{b_s(t)} w_j(s) \right)_{j}$. Therefore, it follows by H\"older's inequality
\begin{align*}
|\langle y, w\rangle| &= \left|\int_0^\infty \langle x(t),\, F(t) \rangle_E \, d\nu(t)\right|\\
&\le \int_0^\infty  \|x(t)\|_E \|F(t)\|_{E^*} \, d\nu(t)\\
&\le \|x\|_{L^p((0,\infty), d\nu; E)} \|F\|_{L^{p'}((0,\infty),d\nu;E^*)}.
\end{align*}

We complete the proof by showing $\|F\|_{L^{p'}((0,\infty),d\nu;E^*)} \le M_1 \|w\|_{E^*}$. Since $1<q'\le p'$ and $p'\ge 2$, applying Minkowski's integral inequality twice, one gets
\begin{align*}
    \|F\|_{L_t^{p'} \left( \ell_n^{q'} \left(L_s^2 \right)  \right)} &= \left\|  \left( \| F_j(t,\cdot) \|_{L_s^2} \right)_j  \right\|_{L^{p'}_t (\ell_n^{q'})}\le \left\|  \left( \| F_j(t,\cdot) \|_{L_s^2} \right)_j  \right\|_{\ell_n^{q'} (L^{p'}_t) }\\
    &= \left( \sum_{j=1}^n \|F_j\|_{L_t^{p'}(L_s^2)}^{q'} \right)^{1/q'}\le \left( \sum_{j=1}^n \|F_j\|_{L_s^2(L_t^{p'})}^{q'} \right)^{1/q'}  \\
    &= \left( \sum_{j=1}^n  \left( \int_I \| F_j(\cdot,s) \|_{L_t^{p'}}^2 \, d\rho(s)   \right)^{q'/2}     \right)^{1/q'}.
\end{align*}
By assumption, $\| F_j(\cdot,s) \|_{L_t^{p'}} = |w_j(s)|\|b_s\|_{L^{p'}(d\nu)}\le M_1 |w_j(s)|$. Insert this into the above estimates:
\begin{align*}
    \|F\|_{L_t^{p'}(E^*)} \le M_1  \left( \sum_{j=1}^n  \left( \int_I |w_j(s)|^2 \, d\rho(s)   \right)^{q'/2}     \right)^{1/q'} = M_1 \|w\|_{E^*},
\end{align*}
completing the proof of Lemma~\ref{lem:diagonal-time}.

\end{proof}

\section{\texorpdfstring{$L^p$}{Lp} holomorphic functional calculus}\label{sec:calculus}

\subsection{Imaginary powers and the boundary calculus}\label{subsec:upper}

In this section, we prove the imaginary-power and multiplier
estimates in Theorem~\ref{thm:A}. We first record the classical
estimates and their consequences that will be used in the proof.
The spectral theorem gives $P_tf\to\Pi_0f$ in $L^2$ as
$t\to\infty$. On finite-coordinate $L^2$ functions, mixed-space
contractivity and Fatou's lemma imply
$\|\Pi_0f\|_{L^r(X;\ell^q)}\le\|f\|_{L^r(X;\ell^q)}$.
Density supplies the consistent contractive projections on all
the mixed spaces below, including the scalar case.

\begin{itemize}
\item \emph{Scalar imaginary powers.}
Recall that $\Pi_0=\mathbf1_{\{0\}}(\mathcal{L})$ is the spectral
projection onto the kernel of $\mathcal{L}$. The scalar
imaginary-power estimate of Carbonaro--Dragi\v cevi\'c \cite[Proposition~11]{CD} states
\begin{equation}\label{eq:scalar-input}
 \|\mathcal{L}^{i\tau}(I-\Pi_0)\|_{L^r(X)\to L^r(X)}
 \le C_r(1+|\tau|)^{1/2}e^{\sigma_r|\tau|},
 \qquad 1<r<\infty,\quad \tau\in\mathbb R.
\end{equation}

\item \emph{The vector-valued calculus.}
For every $1<r,q<\infty$, the injective part of $\mathcal{L}$ on
$L^r(X;\ell^q)$ admits a bounded $H^\infty$ calculus
on some sector of angle strictly less than $\pi/2$;
see \cite[Theorem~4]{Xu}. Its hypotheses hold because the scalar
semigroup is contractively regular by its positive modulus,
analytic by \cite[Corollary~6.1]{HaaseForms}, and $\ell^q$
is a UMD Banach lattice. The initial angle and bound can be
chosen in terms of $r,q$ alone. We use the existence of this
initial calculus and its square-function estimates.

Fix an integer $M\ge2$ and put $\psi(z)=z^Me^{-z}$. For $f=(f_j)_j\in L^r(X;\ell^q)$, consider
\[
 \mathcal{S}f=(\psi(s\mathcal{L})f)_{s>0},
 \qquad
 (\mathcal{S}f)_{j,s}(x)=(\psi(s\mathcal{L})f_j)(x).
\]
Here $L^2(ds/s)$ denotes $L^2((0,\infty),ds/s)$.
The symbol $\psi$ vanishes to positive order at zero and decays
exponentially on sectors strictly smaller than $\pi/2$.
The square-function formula \cite[Proposition~9, formula~(7)]{Xu}
therefore applies as follows:
\begin{equation}\label{um:analysis}
 \|\mathcal{S}f\|_{L^r(X;\ell^q(L^2(ds/s)))}
 \le C_{r,q}\|f\|_{L^r(X;\ell^q)}.
\end{equation}
Also, for a family $v=(v_{j,s}(x))_{j,s}$, write $v_s(x)=(v_{j,s}(x))_j$. Define 
\[
 \mathcal Rv
 =\int_0^\infty\psi(s\mathcal{L})v_s\,\frac{ds}{s}.
\]
It follows from \eqref{um:analysis} at the dual
exponents $(r',q')$ that
\begin{equation}\label{um:synthesis}
 \|\mathcal Rv\|_{L^r(X;\ell^q)}
 \le C_{r,q}
 \|v\|_{L^r(X;\ell^q(L^2(ds/s)))}.
\end{equation}
Indeed, since
\[
 \langle\mathcal Rv,g\rangle
 =\int_0^\infty
   \langle v_s,\psi(s\mathcal{L})g\rangle\,\frac{ds}{s},
\]
using Cauchy--Schwarz in the scale variable, followed
by H\"older's inequality in the sequence and spatial
variables, shows
\[
 |\langle\mathcal Rv,g\rangle|
 \le
 \|v\|_{L^r(X;\ell^q(L^2(ds/s)))}
 \|\mathcal{S}g\|_{L^{r'}(X;\ell^{q'}(L^2(ds/s)))}.
\]

For a compact interval $I\subset(0,\infty)$, let
$\mathcal{S}_I$ denote restriction of $\mathcal{S}$
to scales $s\in I$, and let
\[
 \mathcal R_Iv=\int_I\psi(s\mathcal{L})v_s\,\frac{ds}{s}.
\]
Restriction in the scale variable is contractive,
and $\mathcal R_I$ is obtained by extending its input
by zero outside $I$ before applying $\mathcal R$.
Thus \eqref{um:analysis} and \eqref{um:synthesis}
hold for these restricted operators with constants
independent of $I$.
The estimates also hold for finitely many sequence
coordinates, with constants independent of their number.

\item \emph{Weak estimates.}
The strong estimates above imply the weak-space
estimate needed below. More precisely, for each
$1<r<\infty$, choose two fixed exponents
$r_0,r_1$ with $1<r_0<r<r_1<\infty$.
Real interpolation of \eqref{um:synthesis}
at $r_0$ and $r_1$, with interpolation index
$\infty$, provides
\begin{equation}\label{um:weak-synthesis}
 \|\mathcal Rv\|_{L^{r,\infty}(X;\ell^q)}
 \le C_{r,q}
 \|v\|_{L^{r,\infty}(X;\ell^q(L^2(ds/s)))}.
\end{equation}
The same bound holds for $\mathcal R_I$, uniformly
in $I$ and in the number of sequence coordinates.
\end{itemize}

\begin{proof}[Proof of Theorem~\ref{thm:A}, except sharpness]
We first prove the estimates for $1<p<\min(q,2)$.
Fix $n\ge1$ and a compact interval
$I\subset(0,\infty)$. Recall the notation of
Proposition~\ref{prop:common-contact}:
\[
 \mathcal{H}=L^2(I,ds/s),\qquad E=\ell_n^q(\mathcal{H}),\qquad
 (Kv)_{j,s}=s^{-1/2}v_{j,s},\qquad
 a=\frac1{2p}.
\]
Initially, let
$f\in L^p(X;E)\cap L^q(X;E)\cap L^2(X;E)$.
For every $h>0$, the obstacle decomposition provides
\[
 f=(\mathcal{L}^{1/2}+K)u+g,\qquad
 \|g(x)\|_E\le h,\qquad
 g=h\frac{u}{\|u\|_E}\quad\hbox{on }\Omega=\{u\ne0\}.
\]
Put $Z=\Omega^c$ and $V=K^{1/p'}u$.
We use the coordinate notation
\[
 \begin{aligned}
 u_{j,s}(x)&=u_j(x,s),\\
 V_{j,s}(x)&=V_j(x,s)
 =s^{-1/(2p')}u_{j,s}(x)\\
 &=s^{-(1/2-a)}u_{j,s}(x).
 \end{aligned}
\]
We will use the estimates
\[
 \|g\|_{L^p(X;E)}\le\|f\|_{L^p(X;E)},\qquad
 \mu(\Omega)\le h^{-p}\|f\|_{L^p(X;E)}^p,
\]
\[
 \|g\|_{L^q(X;E)}^q
 \le h^{q-p}\|f\|_{L^p(X;E)}^p,
 \qquad
 \mathcal E(V)\le
 \frac{(p-1)^{p-1}}{p^p}\|f\|_{L^p(X;E)}^p,
\]
which follow from
\eqref{um:obstacle-budget} and
\eqref{eq:energy}.

\par\Needspace{5\baselineskip}\medskip
\noindent\emph{Step 1. The outgoing estimate.}
For the moment, let $(H_s)_{s\in I}$ satisfy the
hypotheses of Lemma~\ref{lem:mixed-mellin}, with
bound $M_0$ in \eqref{um:mellin-size}.
Its reconstruction formula reads
\[
 \mathcal{L}^aH_s(\mathcal{L})V_{j,s}
 =
 \sum_{\theta=\pm\phi_{p,q}}
 \int_0^\infty k_s^\theta(t)
       (e^{-te^{i\theta}\mathcal{L}}-I)V_{j,s}\,dt,
\]
where
\[
 \sup_{s\in I}
 \|t^{1+a}k_s^\theta(t)\|_{L^{p'}(dt/t)}
 \le C_{p,q}M_0.
\]
These identities are first understood in the finite-coordinate Hilbert spectral realization $L^2(X;\ell_n^2(\mathcal{H}))$. The simultaneous obstacle construction gives $V\in\operatorname{Dom}(\mathcal{L}^{1/2})$ in that space, and hence $V\in\operatorname{Dom}(\mathcal{L}^a)$ because $a<1/2$. Since $V_{j,s}=0$ on $Z$, 
\[
 \mathbf1_Z\mathcal{L}^aH_s(\mathcal{L})V_{j,s}
 =
 \sum_{\theta=\pm\phi_{p,q}}
 \int_0^\infty k_s^\theta(t)
       \mathbf1_Ze^{-te^{i\theta}\mathcal{L}}V_{j,s}\,dt.
\]
For each ray and a.e. $y\in Z$, apply Lemma~\ref{lem:diagonal-time} with
\[
 b_s(t)=t^{1+a}k_s^\theta(t),\qquad
 x_{j,s}(t)
 =t^{-a}(e^{-te^{i\theta}\mathcal{L}}V_{j,s})(y).
\]
The required membership of $x(t)$ in
$L^p(dt/t;E)$ follows from \eqref{eq:contact}, since $ap=1/2$. Consequently, after summing the two rays, raising
to the $p$-th power and integrating over $y\in Z$,
\begin{align*}
 &\left\|
 \mathbf1_Z(\mathcal{L}^aH_s(\mathcal{L})V_{j,s})_{j,s}
 \right\|_{L^p(X;E)}^p\\
 &\qquad\le
 C_{p,q}M_0^p
 \sum_{\theta=\pm\phi_{p,q}}
 \int_0^\infty
 \|\mathbf1_Ze^{-te^{i\theta}\mathcal{L}}V\|_{L^p(X;E)}^p
 \,\frac{dt}{t^{1+ap}}.
\end{align*}
By \eqref{eq:contact},
\begin{equation}\label{eq:O}
 \left\|
 \mathbf1_Z(\mathcal{L}^aH_s(\mathcal{L})V_{j,s})_{j,s}
 \right\|_{L^p(X;E)}
 \le C_{p,q}M_0\mathcal E(V)^{1/p}.
\end{equation}

\par\Needspace{5\baselineskip}\medskip
\noindent\emph{Step 2. Weak imaginary powers.}
For $\tau\in\mathbb R$, define
\[
 F_s(\lambda)=\lambda^{i\tau}\psi(s\lambda),
 \qquad \lambda>0,\qquad F_s(0)=0,
\]
and
\[
 (\mathcal T_\tau f)_{j,s}=F_s(\mathcal{L})f_{j,s}.
\]
Thus $\mathcal T_\tau$ applies $\mathcal{L}^{i\tau}\psi(s\mathcal{L})$ separately to each $f_{j,s}$. We claim that
\begin{equation}\label{eq:main-claim}
 \|\mathcal T_\tau f\|_{L^{p,\infty}(X;E)}
 \le C_{p,q}(1+|\tau|)
       e^{\sigma_{p,q}|\tau|}
       \|f\|_{L^p(X;E)}.
\end{equation}
We first establish the estimate at exponent $q$:
\begin{equation}\label{eq:main-claim-1}
 \|\mathcal T_\tau v\|_{L^q(X;E)}
 \le C_q(1+|\tau|)^{1/2}e^{\sigma_q|\tau|}
       \|v\|_{L^q(X;E)}.
\end{equation}

Since $L^2((0,\infty),ds/s)$ is a UMD lattice,
\cite[Corollary~3]{Xu} applies.
Extending $v$ by zero outside $I$, we obtain, for some
$0<\vartheta<\pi/2$,
\begin{equation}\label{eq:Xu-Max}
\left\|
 \left(
 \int_I
 \sup_{\zeta\in\Sigma_\vartheta}
 |P_\zeta v_s|^2\,\frac{ds}{s}
 \right)^{1/2}
 \right\|_{L^q(X)}
 \le C_q\|v\|_{L^q(X;\mathcal{H})}.
\end{equation}
The supremum is taken pointwise in $(x,s)$, and the
constant is independent of $I$.

Now, choose $0<\rho<\sin\vartheta$.
For each $s>0$, the circle $|\zeta-s|=\rho s$ lies in
$\Sigma_\vartheta$. Since
$\psi(s\mathcal{L})=(-1)^Ms^M\partial_s^MP_s$, Cauchy's formula guarantees
\[
 \psi(s\mathcal{L})v_s
 =
 \frac{(-1)^MM!s^M}{2\pi i}
 \int_{|\zeta-s|=\rho s}
 \frac{P_\zeta v_s}{(\zeta-s)^{M+1}}\,d\zeta.
\]
Consequently, $|\psi(s\mathcal{L})v_s(x)|\le C_{M,\rho} \sup_{\zeta\in\Sigma_\vartheta}|P_\zeta v_s(x)|$. Taking the $L^q(X;\mathcal{H})$ norm and applying \eqref{eq:Xu-Max} yields
\begin{equation}\label{eq:square}
 \|(\psi(s\mathcal{L})v_s)_{s\in I}\|_{L^q(X;\mathcal{H})}
 \le C_q\|v\|_{L^q(X;\mathcal{H})}.
\end{equation}
The calculation is initially made for finite-rank simple functions in the scale variable and the estimate extends to every $v\in L^q(X;\mathcal{H})$ by density.

Apply Lemma~\ref{lem:gaussian-amplification} to the single
scalar operator $T=\mathcal{L}^{i\tau}(I-\Pi_0)$ on $L^q(X)$.
Estimate \eqref{eq:scalar-input} then yields
\begin{equation}\label{eq:amplify-L}
\|\mathcal{L}^{i\tau}(I-\Pi_0)\|_{L^q(X;\mathcal{H})\to L^q(X;\mathcal{H})}
 \le C_q(1+|\tau|)^{1/2}e^{\sigma_q|\tau|}.
\end{equation}
Factoring $(\mathcal T_\tau v)_{j,s} = \mathcal{L}^{i \tau} (I-\Pi_0) \psi(s\mathcal{L})v_{j,s}$, using \eqref{eq:amplify-L} and then \eqref{eq:square}, proves for each $j$
\begin{equation}
    \|(\mathcal T_\tau v)_j\|_{L^q(X;\mathcal{H})} \le C_q (1+|\tau|)^{1/2} e^{\sigma_q|\tau|} \|v_j\|_{L^q(X;\mathcal{H})}.
\end{equation}
Finally, since $\|v\|_{L^q(X;E)}^q =\sum_{j=1}^n\|v_j\|_{L^q(X;\mathcal{H})}^q$, one gets
\begin{align*}
    \| \mathcal T_\tau v \|_{L^q(X;E)}^q &= \sum_{j=1}^n \|(\mathcal T_\tau v)_j\|_{L^q(X;\mathcal{H})}^q\\
    &\le C_q (1+|\tau|)^{q/2} e^{q\sigma_q|\tau|} \sum_{j=1}^n\|v_j\|_{L^q(X;\mathcal{H})}^q\\
    &= C_q (1+|\tau|)^{q/2} e^{q\sigma_q|\tau|} \|v\|_{L^q(X;E)}^q
\end{align*}
as desired. This proves \eqref{eq:main-claim-1}. Note that all constants are independent of $n$ and $I$.

We now return to the obstacle decomposition of $f$.
For $\lambda>0$, put
\[
 R_s(\lambda)
 =F_s(\lambda)\bigl[1+(s\lambda)^{-1/2}\bigr],
 \qquad
 H_s(\lambda)=(s\lambda)^{1/2-a}R_s(\lambda),
\]
and assign both symbols the value zero at $\lambda=0$.
The scalar identities
\begin{equation*}
\begin{aligned}
 R_s(\lambda)\lambda^{1/2}
 &=F_s(\lambda)(\lambda^{1/2}+s^{-1/2}),\\
 \lambda^aH_s(\lambda)s^{-(1/2-a)}
 &=\lambda^{1/2}R_s(\lambda)
\end{aligned}
\end{equation*}
give
\begin{equation}\label{um:local-decomp}
 \begin{aligned}
 \mathcal T_\tau f
 &=(R_s(\mathcal{L})\mathcal{L}^{1/2}u_{j,s})_{j,s}
   +\mathcal T_\tau g,\\
 R_s(\mathcal{L})\mathcal{L}^{1/2}u_{j,s}
 &=\mathcal{L}^aH_s(\mathcal{L})V_{j,s}.
 \end{aligned}
\end{equation}
Here the first identity uses
$f_{j,s}=(\mathcal{L}^{1/2}+s^{-1/2})u_{j,s}+g_{j,s}$,
and the second uses
$V_{j,s}=s^{-(1/2-a)}u_{j,s}$. Indeed, $R_s(\lambda) = \lambda^{i\tau}e^{-s\lambda}\bigl[(s\lambda)^M+(s\lambda)^{M-1/2}\bigr]$ is bounded on the positive axis uniformly in $s>0$.
The estimate \eqref{eq:O} will provide the required mixed-norm bound for its action on $Z$.

We next check the hypotheses of Lemma~\ref{lem:mixed-mellin} for this particular family $H_s$. Observe that
\[
 H_s(\lambda)
 =
 \lambda^{i\tau}e^{-s\lambda}
 \bigl[(s\lambda)^{M+1/2-a}+(s\lambda)^{M-a}\bigr].
\]
Both exponents $M+1/2-a$ and $M-a$ are positive.
Thus $H_s$ is continuous on $(0,\infty)$ and
$H_s(e^\cdot)\in L^1(\mathbb R)$.
For either exponent $b=M+1/2-a$ or $b=M-a$,
the substitution $u=s\lambda$ verifies
\[
 \int_0^\infty
 \lambda^{i(\tau-\xi)}e^{-s\lambda}(s\lambda)^b
 \,\frac{d\lambda}{\lambda}
 =
 s^{i(\xi-\tau)}
 \Gamma\bigl(b+i(\tau-\xi)\bigr).
\]
Consequently,
\[
 \begin{aligned}
 &\widehat{H_s(e^\cdot)}(\xi)\\
 &\quad=
 s^{i(\xi-\tau)}
 \left\{
 \begin{aligned}
 &\Gamma(M+1/2-a+i(\tau-\xi))\\
 &\quad+\Gamma(M-a+i(\tau-\xi))
 \end{aligned}
 \right\}.
 \end{aligned}
\]
Since $|s^{i(\xi-\tau)}|=1$, Stirling's formula implies
\[
 |\widehat{H_s(e^\cdot)}(\xi)|
 \le
 C_{p,M}(1+|\xi-\tau|)^{M-a}
 e^{-\pi|\xi-\tau|/2}.
\]
Moreover, $1+|\xi|\le(1+|\tau|)(1+|\xi-\tau|)$ and $|\xi|\le|\tau|+|\xi-\tau|$. Therefore
\begin{align*}
 &(1+|\xi|)^{a+1/2}e^{\sigma_{p,q}|\xi|}
 |\widehat{H_s(e^\cdot)}(\xi)|\\
 &\quad\le
 C_{p,q}(1+|\tau|)^{a+1/2}
 e^{\sigma_{p,q}|\tau|}\\
 &\qquad\times
 (1+|\xi-\tau|)^{M+1/2}
 e^{-(\pi/2-\sigma_{p,q})|\xi-\tau|}.
\end{align*}
The last factor belongs to $L^p(d\xi)$ because
$\sigma_{p,q}<\pi/2$. Taking its norm yields
\begin{equation}\label{um:gamma-bound}
 \sup_{s>0}
 \left\|
 (1+|\xi|)^{a+1/2}e^{\sigma_{p,q}|\xi|}
 \widehat{H_s(e^\cdot)}(\xi)
 \right\|_{L^p(d\xi)}
 \le
 C_{p,q}(1+|\tau|)^{a+1/2}
 e^{\sigma_{p,q}|\tau|}.
\end{equation}

Set $D_\tau=(1+|\tau|)e^{\sigma_{p,q}|\tau|}$. Since $a+1/2<1$ and $\sigma_q\le\sigma_{p,q}$,
$D_\tau$ bounds both the size in
\eqref{um:gamma-bound} and the operator norm in
\eqref{eq:main-claim-1}, up to fixed constants.
Applying \eqref{eq:O} to this family and using
\eqref{um:local-decomp} and \eqref{eq:energy},
one obtains
\begin{equation*}
    \|\mathbf1_Z(R_s(\mathcal{L})\mathcal{L}^{1/2}u_{j,s})_{j,s}\|_{L^p(X;E)}
 \le
 C_{p,q}D_\tau\mathcal E(V)^{1/p}
 \le
 C_{p,q}D_\tau\|f\|_{L^p(X;E)}.
\end{equation*}
Consequently, for every $\lambda>0$, 
\begin{align*}
\mu\{\|\mathcal T_\tau f\|_E>2\lambda\}
 &\le\mu(\Omega)\\
 &\quad+\mu\{\|\mathbf1_Z
 (R_s(\mathcal{L})\mathcal{L}^{1/2}u_{j,s})_{j,s}\|_E>\lambda\}\\
 &\quad+\mu\{\|\mathcal T_\tau g\|_E>\lambda\}.
\end{align*}
The first term on the right is bounded by the obstacle estimate. For the second, use Chebyshev's inequality at
exponent $p$. For the third, use Chebyshev at exponent $q$,
\eqref{eq:main-claim-1}, and the bound for $\|g\|_q$.
This yields
\[
 \mu\{\|\mathcal T_\tau f\|_E>2\lambda\}
 \le C_{p,q}
 \left[
 h^{-p}
 +\frac{D_\tau^p}{\lambda^p}
 +\frac{D_\tau^qh^{q-p}}{\lambda^q}
 \right]\|f\|_{L^p(X;E)}^p.
\]
Choose $h=\lambda/D_\tau$. Then $h^{-p}=\frac{D_\tau^p}{\lambda^p}$, $\frac{D_\tau^qh^{q-p}}{\lambda^q}
 =\frac{D_\tau^p}{\lambda^p}$. So,
\begin{equation*}
\mu\{\|\mathcal T_\tau f\|_E>2\lambda\} \le C_{p,q}\frac{D_\tau^p}{\lambda^p} \|f\|_{L^p(X;E)}^p.
\end{equation*}
This proves \eqref{eq:main-claim} on $L^p(X;E)\cap L^q(X;E)\cap L^2(X;E)$. Density extends this estimate to every $f\in L^p(X;E)$.

We next remove the scale variable $s$. Apply
\eqref{eq:main-claim} to $\mathcal{S}_If_0$, and then
apply \eqref{um:weak-synthesis}. Together with
\eqref{um:analysis}, this implies
\[
 \|\mathcal R_I\mathcal T_\tau\mathcal{S}_If_0
   \|_{L^{p,\infty}(X;\ell_n^q)}
 \le
 C_{p,q}D_\tau\|f_0\|_{L^p(X;\ell_n^q)}.
\]
On the spectral core,
\begin{align*}
 \mathcal R_I\mathcal T_\tau\mathcal{S}_I f_0
 &=
 \int_I
 \psi(s\mathcal{L})\mathcal{L}^{i\tau}\psi(s\mathcal{L})^2f_0\,\frac{ds}{s}\\
 &=
 \mathcal{L}^{i\tau}
 \left[\int_I\psi(s\mathcal{L})^3\,\frac{ds}{s}\right]f_0.
\end{align*}
For every $\lambda>0$,
\[
 \int_0^\infty\psi(s\lambda)^3\,\frac{ds}{s}
 =
 \int_0^\infty u^{3M-1}e^{-3u}\,du
 =
 \frac{\Gamma(3M)}{3^{3M}}
 =:c_M,
\]
whereas the integral is zero at $\lambda=0$.
Exhausting $(0,\infty)$ by compact intervals
therefore shows that $\mathcal R_I\mathcal T_\tau\mathcal{S}_I f_0 \to c_M\mathcal{L}^{i\tau}(I-\Pi_0)f_0$ strongly in $L^2(X;\ell_n^2)$. Weak-type Fatou establishes
\begin{equation}\label{um:weak-powers}
 \begin{gathered}
 \|\mathcal{L}^{i\tau}(I-\Pi_0)\|_
 {L^p(X;\ell_n^q)\to L^{p,\infty}(X;\ell_n^q)}
 \le C_{p,q}(1+|\tau|)e^{\sigma_{p,q}|\tau|},\\
 1<p<\min(q,2).
 \end{gathered}
\end{equation}
Increasing coordinate projections and Fatou's lemma give the same estimate on $\ell^q$.
The constants are locally uniform in $p$ in the initial range:
the contact and Mellin constants are locally bounded, while
the analysis and synthesis bounds are controlled by interpolation
between fixed surrounding exponents.

\par\Needspace{5\baselineskip}\medskip
\noindent\emph{Step 3. Strong imaginary powers in the initial range.}
Put $T_\tau=\mathcal{L}^{i\tau}(I-\Pi_0)$ and fix
$1<p<\min(q,2)$. Choose $\delta_0>0$ sufficiently small that
$1/(1/p\pm\delta_0)$ remain in this range, and set
\[
 \delta=\min\{\delta_0,(1+|\tau|)^{-1}\},
 \qquad
 p_-=\frac1{1/p+\delta},
 \qquad
 p_+=\frac1{1/p-\delta}.
\]
Then $p_-<p<p_+$ and
$1/p=(1/p_-+1/p_+)/2$.
Consequently, \eqref{um:weak-powers} supplies
\[
 A_\pm:=
 \|T_\tau\|_{L^{p_\pm}(X;\ell^q)
                 \to L^{p_\pm,\infty}(X;\ell^q)}
 \le C_{p,q}(1+|\tau|)
                e^{\sigma_{p_\pm,q}|\tau|}.
\]
By the Marcinkiewicz interpolation theorem,
\begin{align*}
\|T_\tau\|_{L^p(X;\ell^q)\to L^p(X;\ell^q)}
 &\le
 2\left[
 p\left(\frac1{p-p_-}+\frac1{p_+-p}\right)
 \right]^{1/p} A_-^{1/2}A_+^{1/2}\\
 &\le C_p\delta^{-1/p}A_-^{1/2}A_+^{1/2}.
\end{align*}
Here the last inequality follows from $(p-p_-)^{-1} + (p_+-p)^{-1} = 2 (p^2 \delta)^{-1}$.

The function
$x\mapsto\max\{\arcsin|1-2x|,\sigma_q\}$
is locally Lipschitz on $(0,1)$. Indeed, away from $x=1/2$
the derivative of $\arcsin|1-2x|$ has absolute value
$1/\sqrt{x(1-x)}$, and taking the maximum with a constant
preserves local Lipschitz continuity. Hence
\[
 \frac{\sigma_{p_-,q}+\sigma_{p_+,q}}2
 \le \sigma_{p,q}+C_{p,q}\delta.
\]
Since $|\tau|\delta\le1$ and
$\delta^{-1/p}\le C_{p,q}(1+|\tau|)^{1/p}$, the preceding
estimates imply
\begin{equation}\label{um:strong-initial}
 \begin{aligned}
 \|T_\tau\|_{L^p(X;\ell^q)\to L^p(X;\ell^q)}
 &\le
 C_{p,q}\delta^{-1/p}(1+|\tau|)
 e^{\sigma_{p,q}|\tau|+C_{p,q}\delta|\tau|}\\
 &\le
 C_{p,q}(1+|\tau|)^{1+1/p}e^{\sigma_{p,q}|\tau|}\\
 &\le
 C_{p,q}(1+|\tau|)^2e^{\sigma_{p,q}|\tau|}.
 \end{aligned}
\end{equation}

\par\Needspace{5\baselineskip}\medskip
\noindent\emph{Step 4. All mixed exponents.}
Suppose first that $q>2$. Estimate \eqref{um:strong-initial}
already covers $1<p<2$. At $p=q$, Fubini and
\eqref{eq:scalar-input} give the scalar bound with angle
$\sigma_q$. Choose $q'<p_0<2$; then
$2/p_0-1<1-2/q$, so $\sigma_{p_0,q}=\sigma_q$.
The mixed-norm Riesz--Thorin theorem \cite[Theorem~2, p.~316]{BP1961} applied in the outer exponent between $p_0$ and $q$
therefore proves the degree-$2$ bound with angle
$\sigma_q=\sigma_{p,q}$ for $2\le p\le q$.
For $p>q>2$, Lemma~\ref{lem:gaussian-amplification} applied to the scalar $L^p$ estimate supplies its bound on $L^p(X;\ell^2)$, while Fubini recovers the same
bound on $L^p(X;\ell^p)$. Applying the same interpolation theorem
with the outer exponent fixed and
$1/q=(1-\vartheta)/2+\vartheta/p$, where $0<\vartheta<1$,
yields the estimate on $L^p(X;\ell^q)$ with angle
$\sigma_p=\sigma_{p,q}$. When $q=2$, the same lemma applies directly for every $p$. Finally, if $1<q<2$, apply the
already proved estimate to $(p',q')$ and $-\tau$, and dualize.
Note that $T_\tau^*=T_{-\tau}$ and $\sigma_{p'}=\sigma_p$, $\sigma_{q'}=\sigma_q$. All these estimates are uniform in the number of coordinates. We have therefore proved
\begin{equation}\label{um:strong-powers}
 \|\mathcal{L}^{i\tau}(I-\Pi_0)\|_{L^p(X;\ell^q)\to L^p(X;\ell^q)}
 \le C_{p,q}(1+|\tau|)^2e^{\sigma_{p,q}|\tau|},
 \qquad 1<p,q<\infty.
\end{equation}

\par\Needspace{5\baselineskip}\medskip
\noindent\emph{Step 5. The boundary multiplier calculus.}
Suppose $\sigma_{p,q}>0$. Recall that $\|\Pi_0\|_{L^p(X;\ell^q)\to L^p(X;\ell^q)}\le1$. The range--kernel decomposition on this reflexive space is
\[
 (I-\Pi_0)L^p(X;\ell^q)
 =\overline{\operatorname{Ran}\mathcal{L}}.
\]
The restriction of $\mathcal{L}$ to this space, with domain
$\operatorname{Dom}(\mathcal{L})\cap\overline{\operatorname{Ran}\mathcal{L}}$,
is densely defined, injective, and has dense range.
It inherits the initial bounded holomorphic calculus recorded
before \eqref{um:analysis}; the restriction and its properties
are precisely those discussed in \cite[Remark~9.7]{HP}.

Let $s>3$. To apply \cite[Theorem~9.4]{HP}, take
$\omega=\sigma_{p,q}$,
$\widetilde v(t)=C_{p,q}(1+|t|)^2$, and
$v(t)=(1+|t|)^s$. The weight $v$ is strongly admissible
(see \cite[Examples~2.4(1)]{HP}), and
\[
 \frac{\widetilde v(t)}{v(t)}
 =C_{p,q}(1+|t|)^{2-s},
 \qquad
 \int_{\mathbb R}\frac{\widetilde v(t)}{v(t)}\,dt
 =\frac{2C_{p,q}}{s-3}<\infty.
\]
Estimate \eqref{um:strong-powers}, with $\tau=-t$, supplies
the required bound for the negative imaginary powers.
Thus \cite[Theorem~9.4]{HP} supplies the boundary
H\"ormander calculus. By \cite[Proposition~9.1(d) and Remark~9.2]{HP}, its multiplier norm is equivalent to our $\cH_{\sigma_{p,q}}^s$. Writing $m(\mathcal{L})f=m(\mathcal{L})(I-\Pi_0)f+m(0)\Pi_0f$ and using $\|I-\Pi_0\|\le2$, we obtain \eqref{eq:horm-main}. For $\theta>\sigma_{p,q}$, \cite[Proposition~9.1(b)]{HP} ensures
\[
 \cH_{\sigma_{p,q}}^s(m)
 \le C_{s,\theta,p,q}\|m\|_{H^\infty(\Sigma_\theta)}.
\]
This proves the $H^\infty$ calculus on every strictly larger sector. The remaining case $p=q=2$ follows directly from the spectral theorem.
\end{proof}

\subsection{Sharpness of the angle}\label{sec:sharpness}

For the injective realization $\mathcal{L}_0$ of $\mathcal{L}$ on
$(I-\Pi_0)L^p(X;\ell^q)$, define
\[
 \begin{aligned}
 &\omega_{H^\infty}(\mathcal{L};L^p(X;\ell^q))\\
 &\quad:=\inf\left\{\alpha\in(0,\pi):
 \begin{gathered}
 \mathcal{L}_0\text{ has a bounded }\\
 H^\infty(\Sigma_\alpha)\text{ calculus}
 \end{gathered}
 \right\}.
 \end{aligned}
\]
Scalar sharpness of the angle $\sigma_p$ is classical for the
Ornstein--Uhlenbeck operator; see
Garc\'ia-Cuerva--Mauceri--Meda--Sj\"ogren--Torrea
\cite[Theorem~2]{GCMST},
Hebisch--Mauceri--Meda \cite[Theorem~2.2]{HMM},
and the discussion following \cite[Theorem~1]{CD}.

The corresponding dimension-free Ornstein--Uhlenbeck
obstruction was independently obtained by van Neerven \cite{vanNeervenOUSeptember}; finite-dimensional target witnesses are also given in
\cite[Proposition~4.5]{vanNeervenOptimalOU}.
We give a separate realization using finite Hamming cubes. The construction produces one fixed symmetric contraction generator, while a translation embedding transfers the classical scalar obstruction to the fibre exponent.

\begin{theorem}\label{prop:umd-sharp-cubes}
There is a probability space $(X,\mu)$ and a nonnegative
selfadjoint scalar generator $\mathcal{L}$, whose semigroup is positive, conservative, and contractive on every scalar $L^r$, such that, for every $1<p,q<\infty$,
\begin{equation}\label{eq:umd-cube-lower}
 \omega_{H^\infty}\bigl(\mathcal{L};L^p(X;\ell^q)\bigr)
 \ge \max\{\sigma_p,\sigma_q\}.
\end{equation}
Its optimal bounded holomorphic semigroup angle on
$L^p(X;\ell^q)$ is exactly $\phi_{p,q}$.
\end{theorem}

\begin{proof}
\emph{Step 1. The complex contractivity obstruction.}
Let $G_n=(\mathbb Z/2\mathbb Z)^n=\{0,1\}^n$ carry
normalized counting measure $\nu_n$, i.e., for $x\in G_n$, $\nu_n(\{x\}) = 2^{-n}$, and define
\[
 T_jf(x)=f(x+e_j),\qquad
 L_n=\sum_{j=1}^n(I-T_j),
\]
where $e_j$ is the $j$th unit coordinate vector and addition is
modulo $2$.
Thus $L_1^{(j)}:=I-T_j$ is the two-point operator $L_1$
acting in the $j$th variable.

The scalar facts used below are classical; see
Ivanisvili--Nazarov
\cite[Lemma~2.1, Subsections~2.2 and~2.4]{IvanisviliNazarov}.
We recall their precise consequences in our normalization.
Since $T_j^2=I$ and the coordinate translations commute, one has by the exponential formula that $e^{-z L_1^{(j)}} = \alpha_zI+\beta_zT_j$, where $\alpha_z=(1+e^{-2z})/2$ and $\beta_z=(1-e^{-2z})/2$. Immediately,
\[
 e^{-zL_n}
 =\prod_{j=1}^n (\alpha_zI+\beta_zT_j).
\]
In particular, with the notation
\[
 N_r(z)=\|e^{-zL_1}\|_{L^r(G_1)\to L^r(G_1)},
\]
the classical tensor-power argument gives (see \cite[Lemma~2.1]{IvanisviliNazarov})
\begin{equation}\label{eq:tensor}
\|e^{-zL_n}\|_{L^r(G_n)\to L^r(G_n)}=N_r(z)^n.    
\end{equation}
Indeed, applying the two-point estimate in each variable
and summing over the remaining variables gives the upper
bound. For the reverse inequality, choose $g$ with
$\|g\|_{L^r(G_1)}=1$ and
$\|e^{-zL_1}g\|_{L^r(G_1)}=N_r(z)$.
For $f_n(x)=\prod_{j=1}^ng(x_j)$, the product measure yields
\[
 \|f_n\|_{L^r(G_n)}=1,\qquad
 e^{-zL_n}f_n(x)
 =\prod_{j=1}^n(e^{-zL_1}g)(x_j),
\]
and hence $\|e^{-zL_n}f_n\|_{L^r(G_n)}=N_r(z)^n$.

We next verify directly the necessary two-point condition (see \cite[Subsection~2.2, equation~(2.7)]{IvanisviliNazarov}):
\begin{equation}\label{eq:Nrz}
N_r(z)\le1 \implies 1-|e^{-2z}|^2 \ge \frac{|r-2|}{r}|1-e^{-4z}|.
\end{equation}
Suppose that $N_r(z)\le1$. For $\zeta\in\mathbb C$ with
$|\zeta|=1$ and $\varepsilon>0$, put $v_\varepsilon=(1+\varepsilon\zeta,\,1-\varepsilon\zeta)$. The two-point exponential formula indicates
\begin{equation*}
e^{-zL_1}v_\varepsilon
 =
 (1+\varepsilon e^{-2z}\zeta,\,
  1-\varepsilon e^{-2z}\zeta).
\end{equation*}
Contractivity on $L^r(G_1)$ therefore implies
\begin{equation}\label{eq:contract}
\frac{|1+\varepsilon e^{-2z}\zeta|^r +|1-\varepsilon e^{-2z}\zeta|^r}{2}
 \le \frac{|1+\varepsilon\zeta|^r +|1-\varepsilon\zeta|^r}{2}.
\end{equation}
For each fixed $\eta\in\mathbb C$, Taylor expansion of $|1\pm\varepsilon\eta|^r = \bigl(1\pm2\varepsilon\Rea\eta +\varepsilon^2|\eta|^2\bigr)^{r/2}$ at $\varepsilon=0$ yields
\[
 \frac{|1+\varepsilon\eta|^r+|1-\varepsilon\eta|^r}{2}
 =
 1+\frac r2\varepsilon^2
 \bigl(|\eta|^2+(r-2)(\Rea\eta)^2\bigr)
 +o(\varepsilon^2).
\]
Subtracting $1$ on both sides of \eqref{eq:contract}, dividing by $r\varepsilon^2/2$, and letting $\varepsilon\downarrow 0$, we obtain
\[
 |e^{-2z}|^2
 +(r-2)\bigl(\Rea(e^{-2z}\zeta)\bigr)^2
 \le
 1+(r-2)(\Rea\zeta)^2.
\]
Using $2(\Rea\eta)^2=|\eta|^2+\Rea(\eta^2)$, this becomes
\[
 \Rea\bigl((r-2)(e^{-4z}-1)\zeta^2\bigr)
 \le r\bigl(1-|e^{-2z}|^2\bigr).
\]
Note that $\zeta^2$ ranges over the unit circle. Taking the supremum over $|\zeta^2|=1$ gives
\[
 |r-2|\,|1-e^{-4z}|
 \le r\bigl(1-|e^{-2z}|^2\bigr),
\]
which proves the necessary condition \eqref{eq:Nrz}.

Now, fix $\phi_r<\theta<\pi/2$. Consider the ray $z = t e^{i \theta}$ ($t>0$). As $t\to 0$,
\begin{equation*}
1- \left|e^{-2te^{i\theta}} \right|^2
 =1-e^{-4t\cos\theta}
 =4t\cos\theta+O(t^2),
 \qquad
 |1-e^{-4te^{i\theta}}|
 =4t+O(t^2).
\end{equation*}
Since $\cos\theta<\cos\phi_r=\frac{|r-2|}{r}$, the necessary condition in \eqref{eq:Nrz} fails for all sufficiently small $t>0$. Fix one such $t$. It follows that $N_r(te^{i\theta}) > 1$, and \eqref{eq:tensor} guarantees
\begin{equation}\label{eq:umd-scalar-cube-obstruction}
 \sup_{n\ge1}
 \|e^{-te^{i\theta}L_n}\|_{L^r(G_n)\to L^r(G_n)}
 =\infty.
\end{equation}

\medskip
\noindent\emph{Step 2. Construction of the space and generator.}
Put
\[
 X=\bigsqcup_{n\ge1}G_n,\qquad
 \mu|_{G_n}=2^{-n}\nu_n.
\]
Use coordinates $(n,x) \in X$, where $n\ge 1$ and $x\in G_n$. Thus $\mu((n,x)) = 2^{-2n}$ and $\mu(X)=\sum_{n\ge1}2^{-n}=1$. Writing $f_n=f|_{G_n}$,
we have
\[
 \|f\|_{L^r(X)}^r
 =\sum_{n\ge1}2^{-n}\|f_n\|_{L^r(G_n)}^r.
\]
Define $\mathcal{L}$ on $L^2(X)$ by
\[
(\mathcal{L}f)|_{G_n}=L_nf_n,\qquad n\ge 1,
\]
with domain
\[
 \operatorname{Dom} \mathcal{L}
 =
 \left\{
 f\in L^2(X):
 \sum_{n\ge1}2^{-n}\|L_nf_n\|_{L^2(G_n)}^2<\infty
 \right\}.
\]

We next show that $\mathcal{L}$ is nonnegative and selfadjoint. For $S\subseteq\{1,\ldots,n\}$, the Walsh functions
\[
 \chi_S(x)=(-1)^{\sum_{j\in S}x_j}
\]
form an orthonormal basis of $L^2(G_n)$ and satisfy $L_n\chi_S=2|S|\chi_S$. Consequently, $2^{n/2} \mathbf 1_{G_n} \chi_S$ form an orthonormal basis of $L^2(X)$ in which $\mathcal{L}$ is diagonal with nonnegative real entries. With the domain above, $\mathcal{L}$ is therefore
nonnegative and selfadjoint.

For $t\ge0$, $\alpha_t,\beta_t\ge0$ and
$\alpha_t+\beta_t=1$. Each operator $e^{-t L_1^{(j)}}$, i.e.,
$\alpha_tI+\beta_tT_j$, is positive and conservative,
and contracts every scalar $L^r$. The same properties
hold for their products $e^{-tL_n}$ and hence for
\[
(e^{-t\mathcal{L}}f)|_{G_n}=e^{-tL_n}f_n.
\]
Approximation by finitely many blocks gives strong continuity
on $L^r(X)$ for $1\le r<\infty$.

\medskip
\noindent\emph{Step 3. Transferring the fibre exponent.}
For a sequence-valued function $F=(F_k)$ on $G_n$,
we let $e^{-zL_n}$ act coordinatewise:
\[
 (e^{-zL_n}F)_k=e^{-zL_n}F_k.
\]
The main aim of this step is to prove
\begin{equation}\label{eq:umd-cube-mixed-obstruction}
 \|e^{-zL_n}\|_{L^p(G_n;\ell^q)\to L^p(G_n;\ell^q)} \ge
 \max\left\{ N_p(z)^n, N_q(z)^n \right\}.
\end{equation}
Note that the first lower bound follows by testing on
$F(x)=f(x)\mathbf e_1$, where $\mathbf e_1 = (1,0,\dots) \in \ell^q$. 

For the second lower bound, identify the first $|G_n|=2^n$ coordinates
of $\ell^q$ with the elements of $G_n$, in any fixed order, and define
\[
 (\mathcal J_nf)(x)_y
 =|G_n|^{-1/q}f(x+y),\qquad x,y\in G_n,
\]
with all remaining coordinates zero. For each fixed $x$, translation permutes $G_n$ (recall that the addition is modulo 2), so
\begin{align*}
\|\mathcal J_nf(x)\|_{\ell^q}^q
 &=\frac1{|G_n|}\sum_{y\in G_n}|f(x+y)|^q\\
 &=\frac1{|G_n|}\sum_{y\in G_n}|f(y)|^q
 =\|f\|_{L^q(G_n)}^q.
\end{align*}
Since $\nu_n(G_n)=1$, we get $\|\mathcal J_nf\|_{L^p(G_n;\ell^q)} =\|f\|_{L^q(G_n)}$. 

For $y\in G_n$, write $T_yf(x)=f(x+y)$.
A direct calculation shows 
\begin{equation*}
L_nT_yf(x) =\sum_{j=1}^n \bigl(f(x+y)-f(x+e_j+y)\bigr) = T_yL_nf(x).
\end{equation*}
Thus $e^{-zL_n}T_y=T_ye^{-zL_n}$. Applying this identity to each sequence coordinate yields
\begin{align*}
(e^{-zL_n}\mathcal J_nf)(x)_y
 &=|G_n|^{-1/q}e^{-zL_n}(T_yf)(x)\\
 &=|G_n|^{-1/q}(e^{-zL_n}f)(x+y)\\
 &=(\mathcal J_n(e^{-zL_n}f))(x)_y.
\end{align*}
Consequently, $\|e^{-zL_n}\mathcal J_nf\|_{L^p(G_n;\ell^q)}=\|e^{-zL_n}f\|_{L^q(G_n)}$.

It follows that
\begin{equation*}
\frac{\|e^{-zL_n}\mathcal J_nf\|_{L^p(G_n;\ell^q)}}{\|\mathcal J_nf\|_{L^p(G_n;\ell^q)}} = \frac{\|e^{-zL_n}f\|_{L^q(G_n)}}{\|f\|_{L^q(G_n)}}.
\end{equation*}
Taking the supremum over $f\ne0$ yields $\|e^{-zL_n}\|_{L^p(G_n,\ell^q) \to L^p(G_n,\ell^q)} \ge N_q(z)^n$, completing the proof of \eqref{eq:umd-cube-mixed-obstruction}.

\medskip
\noindent\emph{Step 4. The exact semigroup angle.}
Observe that $e^{-z\mathcal{L}}f$ is supported in $G_n$ if $f$ is supported in $G_n$. For a sequence-valued function $F$ with all entries supported in $G_n$,
\begin{equation*}
\begin{aligned}
 \|F\|_{L^p(X;\ell^q)}
 &=2^{-n/p}\|F\|_{L^p(G_n;\ell^q)},\\
 \|e^{-z\mathcal{L}} F\|_{L^p(X;\ell^q)}
 &=2^{-n/p}\|e^{-zL_n} F\|_{L^p(G_n;\ell^q)}.
\end{aligned}
\end{equation*}
Thus any bounded extension of $e^{-z\mathcal{L}}$ to
$L^p(X;\ell^q)$ would bound all the mixed norms
in \eqref{eq:umd-cube-mixed-obstruction} uniformly in $n$.

Suppose that $\phi_{p,q}<\pi/2$ and that the real semigroup
had a bounded holomorphic extension to a sector
$\Sigma_\beta$ with $\beta>\phi_{p,q}$.
Choose $r\in\{p,q\}$ such that $\phi_r=\phi_{p,q}$ and pick
\[
 \phi_r<\theta<\min\{\beta,\pi/2\}.
\]
By \eqref{eq:umd-scalar-cube-obstruction}, there is a fixed
$z=te^{i\theta}$ for which the scalar $L^r(G_n)$ norms
are unbounded. Equation~\eqref{eq:umd-cube-mixed-obstruction} gives the
same conclusion for the mixed norms. This contradicts boundedness of the extension at $z$. Proposition~\ref{prop:umd-complex-contractivity} supplies the reverse bound for the holomorphic angle. Hence the exact angle is $\phi_{p,q}$. 

\medskip
\noindent\emph{Step 5. The kernel projection.}
The Walsh eigenvalue formula shows that $\ker L_n$
consists precisely of the constant functions.
Its orthogonal projection is therefore the normalized
block average. Hence
\[
 (\Pi_0F)(x)=\int_{G_n}F(y)\,d\nu_n(y),
 \qquad x\in G_n.
\]
For sequence-valued functions this formula is understood
coordinatewise. The triangle inequality and H\"older's
inequality give
\[
 \|\Pi_0 F(x)\|_{\ell^q}^p = \left\|\int_{G_n}F\,d\nu_n\right\|_{\ell^q}^p
 \le
 \left(\int_{G_n}\|F\|_{\ell^q}\,d\nu_n\right)^p
 \le
 \int_{G_n}\|F\|_{\ell^q}^p\,d\nu_n.
\]
Multiplying by $2^{-n}$ and summing over $n\ge 1$ proves
\[
 \|\Pi_0F\|_{L^p(X;\ell^q)}
 \le\|F\|_{L^p(X;\ell^q)}.
\]

\medskip
\noindent\emph{Step 6. The lower holomorphic calculus angle.}
If $\sigma_{p,q}=0$, the asserted lower bound is immediate.
Otherwise, suppose that the injective realization of $\mathcal{L}$ had
a bounded $H^\infty(\Sigma_\alpha)$ calculus, with bound
$C$, for some $0<\alpha<\sigma_{p,q}$. Choose $r\in\{p,q\}$ and $\theta$ such that
\[
 \alpha<\sigma_r,\qquad
 \phi_r<\theta<\pi/2-\alpha,
\]
and choose the fixed $z=te^{i\theta}$ supplied by
\eqref{eq:umd-scalar-cube-obstruction}.

For $\delta>0$, define
\[
 m_{z,\delta}(\lambda)
 =\frac{\lambda}{\lambda+\delta}e^{-z\lambda}.
\]
Since $\alpha<\pi/2$, $\Rea \lambda >0$ for every $\lambda\in\Sigma_\alpha$. Moreover, $\theta+\alpha<\pi/2$ implies
\[
 |e^{-z\lambda}|
 \le e^{-c|z||\lambda|},
 \qquad c=\cos(\theta+\alpha)>0.
\]
Consequently, $m_{z,\delta}\in H^\infty_0(\Sigma_\alpha)$ and $\|m_{z,\delta}\|_{H^\infty(\Sigma_\alpha)}\le1$. The assumed calculus therefore yields
\begin{equation*}
\|m_{z,\delta}(\mathcal{L})(I-\Pi_0)F\|_{L^p(X;\ell^q)}
 \le2C\|F\|_{L^p(X;\ell^q)}.
\end{equation*}
Fix $n$ and take $F$ supported in $G_n$, with finitely
many nonzero sequence coordinates.
On each positive eigenvalue of $L_n$,
$m_{z,\delta}(\lambda)\to e^{-z\lambda}$ as
$\delta\to 0$, while $m_{z,\delta}(0)=0$.
There are only finitely many eigenvalues, so passage to the limit gives
\[
 \|e^{-zL_n}(I-\Pi_{0,n})F\|_{L^p(G_n;\ell^q)}
 \le2C\|F\|_{L^p(G_n;\ell^q)},
\]
where
\[
 \Pi_{0,n}F=\int_{G_n}F\,d\nu_n.
\]
Since $e^{-zL_n}\Pi_{0,n}=\Pi_{0,n}$ and
$\|\Pi_{0,n}\|\le1$, we obtain
\[
 \|e^{-zL_n}F\|_{L^p(G_n;\ell^q)}
 \le(2C+1)\|F\|_{L^p(G_n;\ell^q)}.
\]
The constant is independent of $n$ and of the number
of sequence coordinates.
This contradicts
\eqref{eq:umd-scalar-cube-obstruction} and
\eqref{eq:umd-cube-mixed-obstruction}, proving
\eqref{eq:umd-cube-lower}.

\end{proof}

\begin{proof}[Completion of the proof of Theorem~\ref{thm:A}]
The upper bound already proved and
Theorem~\ref{prop:umd-sharp-cubes} imply
\[
 \omega_{H^\infty}\bigl(\mathcal{L};L^p(X;\ell^q)\bigr)
 =\sigma_{p,q}.
\]
The same example also proves sharpness for every finite boundary
Sobolev order. Indeed, for $0<\alpha<\sigma_{p,q}$, choose $r$ and $z$ as in the preceding proof. Since $|\arg z|<\pi/2-\alpha$, the function
$m_z(\lambda)=e^{-z\lambda}$ is bounded on a sector strictly
larger than $\Sigma_\alpha$. Direct estimates on annuli give
$\cH_\alpha^s(m_z)<\infty$ for every finite $s$, whereas $m_z(\mathcal{L})$ is unbounded on the mixed space. Thus no finite increase in the boundary Sobolev order permits a universally smaller angle.
\end{proof}

\section{Applications}\label{sec:applications}

We apply Theorem~\ref{thm:A} to two classical consequences
of holomorphic functional calculus: lattice square functions
and sectorial maximal functions.
The additional information here is the full strict angular
range determined by both mixed exponents.
For the maximal inequality, the universal angular threshold
is optimal.

\subsection{Square functions on complex-time rays}
\label{subsec:rotated-square}

We first obtain square-function estimates uniformly on
every closed subsector of $\Sigma_{\phi_{p,q}}$.

\begin{corollary}\label{cor:rotated-square}
Let $\mathcal{L}$ be as in Theorem~\ref{thm:A}, let
$1<p,q<\infty$, and let $0<\gamma<\phi_{p,q}$.
For every $f=(f_j)_j\in L^p(X;\ell^q)$ and
$|\vartheta|\le\gamma$,
\begin{equation}\label{eq:rotated-square}
 \begin{gathered}
 \left\|
 \left(
 \left[\int_0^\infty
 \left|te^{i\vartheta}\mathcal{L}
       e^{-te^{i\vartheta}\mathcal{L}}f_j\right|^2
 \frac{dt}{t}\right]^{1/2}
 \right)_j
 \right\|_{L^p(X;\ell^q)}
 \\
 \simeq_{p,q,\gamma}
 \|(I-\Pi_0)f\|_{L^p(X;\ell^q)}.
 \end{gathered}
\end{equation}

\end{corollary}

\begin{proof}
Choose $\sigma_{p,q}<\nu_0<\nu_1<\pi/2-\gamma$ and put
$\psi_\vartheta(z)=e^{i\vartheta}ze^{-e^{i\vartheta}z}$.
For $|\vartheta|\le\gamma$ and $z\in\Sigma_{\nu_1}$,
\[
 |\psi_\vartheta(z)|
 \le |z|e^{-\cos(\nu_1+\gamma)|z|}
 \le C_{\nu_1,\gamma}\min\{|z|,|z|^{-1}\}.
\]
Theorem~\ref{thm:A} supplies a bounded calculus on
$\Sigma_{\nu_0}$. The general-symbol square-function
estimate in \cite[Proposition~9, formula~(7) and its proof]{Xu},
extending \cite[Corollary~6.7]{CDMY}, therefore gives
the upper bound uniformly in $\vartheta$.
Its lattice hypotheses hold because $L^p(X;\ell^q)$
has nontrivial convexity and finite concavity, independently
of the number of coordinates.
The kernel contributes nothing since $\psi_\vartheta(0)=0$.

For the reverse estimate, the spectral theorem gives
\[
 \langle(I-\Pi_0)f,g\rangle
 =4\int_0^\infty
 \langle\psi_\vartheta(t\mathcal{L})f,
        \psi_{-\vartheta}(t\mathcal{L})g\rangle
 \,\frac{dt}{t},
\]
since
$\int_0^\infty\psi_\vartheta(t\lambda)^2\,dt/t=1/4$
for $\lambda>0$.
Cauchy--Schwarz in $t$, H\"older's inequality, and the
upper estimate at $(p',q')$ and $-\vartheta$ give the
reverse bound by duality, because
$\phi_{p',q'}=\phi_{p,q}$.
This completes the proof.
\end{proof}

\subsection{Sectorial maximal functions}
\label{subsec:sectorial-maximal}

The following is the mixed-norm counterpart of
\cite[Theorem~12(i)]{CD}, with the angular range supplied
by \eqref{eq:power-main}.

\begin{corollary}\label{cor:sectorial-maximal}
Let $\mathcal{L}$ be as in Theorem~\ref{thm:A}, let
$1<p,q<\infty$, and let $0<\gamma<\phi_{p,q}$.
For $f=(f_j)_j\in L^p(X;\ell^q)$, define the sectorial
maximal function coordinatewise by
\[
 (\mathcal M_\gamma f)_j(x)
 =\sup_{z\in\Sigma_\gamma}|P_zf_j(x)|.
\]
Then
\begin{equation}\label{eq:sectorial-maximal}
 \|\mathcal M_\gamma f\|_{L^p(X;\ell^q)}
 \le C_{p,q,\gamma}\|f\|_{L^p(X;\ell^q)}.
\end{equation}
The universal angular threshold $\phi_{p,q}$ cannot be increased.
\end{corollary}

\begin{proof}
The real semigroup is contractively regular, and $\ell^q$
is a UMD Banach lattice. We therefore use the ergodic
maximal inequality and the lattice Mellin argument in
\cite[Theorem~1 and the proof of Theorem~2]{Xu}.
Writing
\[
 (\mathcal M_{\mathrm{erg}}f)_j
 =\sup_{t>0}\left|\frac1t\int_0^tP_sf_j\,ds\right|,
\]
we have
$\|\mathcal M_{\mathrm{erg}}f\|_{L^p(X;\ell^q)}
 \le C_{p,q}\|f\|_{L^p(X;\ell^q)}$.
The cited Mellin representation, originating in
\cite[Section~3, Theorem~7]{Cowling}, gives
\[
 \mathcal M_\gamma f
 \le \mathcal M_{\mathrm{erg}}f+
 C_\gamma\int_{\mathbb R}
 e^{-(\pi/2-\gamma)|\tau|}
 |\mathcal{L}^{i\tau}(I-\Pi_0)f|\,d\tau.
\]
The ergodic average retains the kernel component, while
the Mellin remainder vanishes at zero.
Minkowski's inequality and \eqref{eq:power-main} imply
\[
 \begin{split}
 \|\mathcal M_\gamma f\|_{L^p(X;\ell^q)}
 &\le C_{p,q,\gamma}
 \left[
 1+\int_{\mathbb R}(1+|\tau|)^2
 e^{-(\phi_{p,q}-\gamma)|\tau|}\,d\tau
 \right]
 \|f\|_{L^p(X;\ell^q)}.
 \end{split}
\]
The integral is finite since $\gamma<\phi_{p,q}$.
The standard maximal extension argument gives the estimate
for every $f\in L^p(X;\ell^q)$.

Finally, the maximal estimate bounds each $P_z$ on its
sector, so Theorem~\ref{prop:umd-sharp-cubes} excludes a
larger universal angular threshold.
For $p=q=2$, every larger sector contains $z$ with
$\Rea z<0$, and the cube eigenvalues $2n$ give
$|e^{-2nz}|\to\infty$.
\end{proof}

\section*{AI Declaration}
During the preparation of this article, the author used ChatGPT (OpenAI; GPT-6) for mathematical discussion, assistance in checking proof arguments, literature searches, and revision of the language and exposition. This assistance included discussion and revision of the initial formulations of the obstacle decomposition and contact inequality arguments, as well as the arguments used in Lemma~2.4 and Proposition~2.7 and Theorem~3.1.
The author takes full responsibility for the mathematical content, proofs, references, and final presentation of the article.

\medskip
\noindent
{\bf Acknowledgments.}
The author would like to thank Professor Adam Sikora for consistent support and helpful discussion. The author is grateful to Professor Jan van Neerven for communicating his unpublished manuscripts, clarifying the chronology of his independent results on vector-valued Ornstein--Uhlenbeck semigroups, and for his generous and constructive correspondence concerning the relationship between the two works. D.He was supported by the China Postdoctoral Science Foundation No.~2026T190831 and the NNSF of China No.~42381947.

\bibliographystyle{abbrv}
\bibliography{references}
\end{document}